\documentclass[12pt, leqno]{article}

\usepackage{wrapfig}
\usepackage{euscript}

\usepackage{tikz, tikz-3dplot, pgfplots}
\usetikzlibrary{decorations.markings}

\usepackage{fullpage}

\usepackage[T1]{fontenc}

\usepackage{amsxtra}
\usepackage{amsmath}
\usepackage{amssymb}
\usepackage{amsfonts}
\usepackage{amsthm}
\usepackage{mathabx}

\usepackage[all,2cell]{xy}
\UseAllTwocells
\usepackage{mathrsfs}
\usepackage{enumitem}
\usepackage{subfigure}
\usepackage{float}
 \usepackage{caption}

\usepackage[latin1]{inputenc}
\usepackage{amsxtra}
\usepackage{amsmath}
\usepackage{amssymb}
\usepackage{amsfonts}
\usepackage[all,2cell]{xy}
\UseAllTwocells
\usepackage{mathrsfs}
\usepackage{amsthm}
\usepackage{enumitem}
\usepackage{subfigure}
  \usepackage{epigraph}

\newtheorem{thm}[equation]{Theorem}
 
\newtheorem{cor}[equation]{Corollary}

\newtheorem{prop}[equation]{Proposition}
\newtheorem{para}[equation]{}
\newtheorem{conj}[equation]{Conjecture}

\newtheoremstyle{example}{\topsep}{\topsep}%
{}
{}
{\bfseries}
{.}
{2pt}
{\thmname{#1}\thmnumber{ #2}\thmnote{ #3}}

\theoremstyle{example}

\newtheorem{Defi}[equation]{Definition}

\newtheorem{rem}[equation]{Remark}
\newtheorem{rems}[equation]{Remarks}

\newtheorem{exas}[equation]{Examples}
\newtheorem{ex}[equation]{Example}

\newtheorem{parab}[equation]{}

\newcounter{sparagraphcounter}[subsection]

   \numberwithin{equation}{subsection}

\def\Ac{\mathcal{A}}

\def\Kc{\mathcal{K}}
\def\Dc{\mathcal{D}}

\def\Fc{\mathcal{F}}
\def\Gc{\mathcal{G}}

\def\Lc{\mathcal{L}}

\def\Oc{\mathcal{O}}
\def\Pc{\mathcal{P}}

\def\Rc{\mathcal{R}}
\def\Sc{\mathcal{S}}

\def\Vc{\mathcal{V}}
\def\Wc{\mathcal{W}}

\def\pen{\mathfrak{p}}

\def\gen{\mathfrak{g}}
\def\men{\mathfrak{m}}

\def\hen{\mathfrak{h}}
\def\ten{\mathfrak{t}}

\def\pen{\mathfrak{p}}

 \def\AAA{\mathbb{A}}
 \def\BB{\mathbb{B}}

\def\CC{\mathbb{C}}

\def\FF{\mathbb{F}}
\def\GG{\mathbb{G}}
\def\KK{\mathbb{K}}

\def\PP{\mathbb{P}}
\def\RR{\mathbb{R}}
\def\SSS{\mathbb{S}}
\def\ZZ{\mathbb{Z}}
\def\QQ{\mathbb{Q}}

 \def\chib{{\boldsymbol{\chi}}}

\def\eps{\varepsilon}

\def\ev{{\bar 0}}
\def\od{{\bar 1}}

\renewcommand\1{{\bf{1}}}

\def\Aut{\text{\rm Aut}}

\def\be{\begin{equation}}
\def\ee{\end{equation}}
\def\beg{\begin{gathered}}
\def\eng{\end{gathered}}

\def\deg{\on{deg}}
\def\Der{\on{Der}}

\def\End{\on{End}}
 \def\Ext{\on{Ext}}

\def\gen{\mathfrak{g}}

\def\FL{\on{FL}}

\def\Hom{\on{Hom}}

\def\Id{\text{\rm Id}}

\def\ISVect{{\on{1-SVect}}}
\def\ISvect{{\on{1-SVect}}}
\def\IAb{{ \on{1-Ab}}}
\def\ISAb {{  \on{1-SAb}}}

\def\Ker{\mathop{\mathrm{Ker}}\nolimits}
\def\k{\mathbf{k}}

\def\Lie{{\on{Lie}}}
\def\lra{\longrightarrow}

\def\ol{\overline}
\def\on{\operatorname}

 \def\SA{{\Sc\Ac}_\k}
 \def\SSch{{\Sc\Sc ch}_\k}
 \def\SVect{\on{SVect}}
 \def\SCat{\on{SCat}}

\def\Spec{{\on{Spec}}}
\def\SLRS{{\Sc\Lc\Rc\Sc}_\k}
\def\SHo{{\on{SHo}}}
\def\Spin{\on{Spin}}
\def\st{{\on{st}}}

\def\Sym{\on{Sym}}

\def\ul{\underline}

\def\Vect{\on{Vect}}
\def\IVect{\on{1-Vect}}

 \def\wt{\widetilde}

\numberwithin{paragraph}{subsection}

\begin{document} 

\title{Supergeometry in mathematics and physics  }
\author{   Mikhail  Kapranov  }
\maketitle


\tableofcontents

\vfill\eject

\numberwithin{equation}{subsection}

\addtocounter{section}{-1}
\section{Introduction} 

 \epigraph{In die Traum- und Zaubersph\"are   \\
Sind wir, scheint es, eingegangen.}
        { Goethe, {\em Faust I }}

 Supergeometry is a geometric tool which is supposed to describe supersymmetry  
 (``symmetry between bosons and fermions'')   in   physics.  In the impressive
 arsenal of physico-geometric tools, it 
   occupies a special,  somewhat mysterious place. 
        
        \vskip .2cm

            To an outside observer, it projects a false sense of accessibility, presenting itself as a simple modification
            of the familiar formalism.   However, this is only an appearance.
            Like a bewitched place, it is protected not by barriers but
            by something  less tangible and therefore much more powerful. 
            It is not even so clear where this place is localized.
            A naive attempt to ``forge ahead''  
           may  encounter little resistance but may also  end up missing the real point.

            \vskip .2cm

        Indeed, the physical origin 
         of supergeometry lies in the comparison of the Bose and Fermi statistics for
        identical particles. So
      mathematically, supergeometry  zooms in  on our  mental habits related to   questions of commutativity and identity. 
        These habits predate quantum physics and therefore are deeply ingrained,
        but they do not correspond to the ultimate physical reality. The new way of thinking, offered by supergeometry
        in mathematics and by supersymmetry in physics,  requires  changes  
        that cannot be reduced to  simple recipes (such as 
          ``introducing $\pm$ signs here and there'' which is, in the popular mind, the best known
          feature of supergeometry).
          
          \vskip .2cm
          
           Put differently, if introducing $\pm$ signs
                   can take us so far, then   something is happening that we truly don't understand. 
                   In addition,  one can almost say that mathematicians and physicists
       mean  different things when they speak about supergeometry. This contributes
       to the enduring feeling of wonder and mystery surrounding this  subject. 
                
       \vskip .2cm

        The present article (chapter) is subdivided into three parts.
         Following the specifications for the  volume, \S 1 presents
        a very short but self-contained exposition of supergeometry as it is  generally understood by mathematicians. 
        Not surprisingly, the presentation revolves around the Koszul sign rule \eqref{eq:supercomm}
        which appears first at the level of elements of a ring, and later again at the level of categories. 
         More systematic expositions can be found in the books \cite{Freed-super, M} and in many articles, esp. 
    \cite{deligne-morgan}.  
    
    \vskip .2cm
    
    In \S 2 we discuss the aspects of supergeometry that are used by physicists in relation to
    supersymmetry. From the mathematical point of view, this amounts to much more than the study of
    super-manifolds or of the Koszul rule. 
      Some of the basic references are \cite{deligne-freed, Freed-super, M, WB}. 
  The entry point for a mathematician here could be found in
     the idea of taking  natural ``square roots'' of familiar mathematical and physical quantities.
     
     In our presentation, we introduce the  abstract concept of a quadratic space (data of $\Gamma$-matrices)
     of which  various situations involving spinors form particular cases. This has the advantage of
     simplifying the general discussion and also of relating the subject to the mathematical theory
     of intersections of quadrics. In particular, the familiar dichotomy of complete intersections of quadrics
     vs. non-complete ones has a direct significance for  understanding the ``constraints" 
     which are usually imposed on super-fields. 
     
          \vskip .2cm
          
        Finally, \S 3 is devoted to   an attempt to uncover deeper roots for the 
        mysterious power of the super-geometric formalism, of which
        the remarkable consistency of its sign rules is just one manifestation. 
      It seems that the right language to speak about such things is given by homotopy theory.
      Indeed, this theory  provides a  systematic modern way to talk about the issues of identity:
      instead of saying that two things are ``the same'', we say that they are ``homotopic'', and specify
      the homotopy (the precise reason why they should be considered the same). We
      can also consider homotopies between homotopies and so on. 
        From this point of view, the group $\{\pm 1\}$ of  signs howering over supergeometry, is nothing
        but $\pi_1^\st$,  the first homotopy group of the {\em sphere spectrum} $\SSS$. As emphasized by
        Grothendieck in a more categorical language,
                $\SSS$ can  be seen as the most fundamental homotopy commutative object. This suggests that 
        ``mining the sphere spectrum'' beyond the first level should lead to generalizations of supergeometry
         (and possibly, supersymmetry) involving not just  two types of quantities (even and odd) 
         subject to $\pm$ sign rules 
          but,  for instance,  24 ``sectors'' (types of higher level categorical objects)
           that would  account for
         $\pi_3^\st=\ZZ/24$ and obey some multi-level sign rules 
         that could use  higher roots of 1 as well.

        \vskip .2cm
        
        My understanding of supergeometry owes a lot to lectures and writings of Y. I.  Manin. 
        In particular, the idea of square roots provided by the  super  formalism,  was learned from him long time ago. 
        The homotopy-theoretic considerations of \S 3 were stimulated by the  joint work with N. Ganter \cite{GK}. 
        I would also like to thank M. Anel,. G. Catren, N. Ganter, N. Gurski and Y. I. Manin for remarks  and suggestions 
        on the preliminary versions of this
        text. 
         This work was supported by World Premier International Research Center Initiative (WPI Initiative), MEXT, Japan.


\section{Supergeometry as understood by mathematicians}

\subsection{Commutative superalgebras} For a mathematician, supergeometry is the study of supermanifolds and superschemes: objects
whose rings of functions are commutative superalgebras. 

\vskip.2cm

We fix a field $\k$ of characteristic $\neq 2$. By an {\em associative superalgebra} over $\k$ one means simply a $\ZZ/2$-graded
associative algebra $A = A^{\ev}\oplus A^\od$. Elements of $A^\ev$ are called {\em even} (or {\em bosonic}), elements of
$A^\od$ are called {\em odd} (or {\em fermionic}). 

An associative superalgebra $A$ is called {\em commutative}, if it satisfies the Koszul sign rule for
commutation of  homogeneous elements:
\be\label{eq:supercomm}
ab \,\,=\,\, (-1)^{\deg(a) \cdot \deg(b)} ba. 
\ee
The terms ``commutative superalgebra'' and ``supercommutative algebra'' are used without distinction. 

Often, one considers $\ZZ$-graded supercommutative algebras $A=\bigoplus_{i\in\ZZ} A^i$, with the same
commutation rule \eqref{eq:supercomm} imposed on homogeneous elements. 

\begin{exas}\label{exas:supercomm}
(a) Any usual commutative algebra $A$ becomes supercommutative, if put in degree $\ev$. The most important example
$A=\k[x_1, \cdots, x_m]$ (the polynomial algebra). 

\vskip .2cm

(b) The exterior (Grassmann) algebra $\Lambda[\xi_1, \cdots \xi_n]$ over $\k$, generated by the symbols $\xi_i$ of degree $\od$,
subject only to the relations
\[
\xi_i^2=0, \quad \xi_i\xi_j=-\xi_j\xi_j, \,\,\, i\neq j, 
\]
is supercommutative. It can be seen as the {\em free supercommutative algebra} on the $\xi_i$: the relations imposed are the minimal ones to
ensure supercommutativity. 

\vskip .2cm

(c) If $A, B$ are the commutative superalgebras, then the tensor product $A\otimes_\k B$ with the standard product grading
$\deg(a\otimes b) =\deg(a) +\deg (b)$ and multiplication
\[
(a_1\otimes b_1) \cdot (a_2\otimes b_2) \,\,\,=\,\,\, (-1)^{\deg(b_1) \cdot \deg(a_2)} (a_1\cdot b_1)\otimes (a_2\cdot  b_2)
\]
($a_i, b_i$ homogeneous), is again a supercommutative algebra. For example,
\[
\k[x_1, \cdots, x_m] \otimes \Lambda[\xi_1, \cdots, \xi_n], \quad \deg(x_i)=\ev, \,\,\,\deg(\xi_j)=\od,
\]
is a commutative superalgebra. This is a general form of a free commutative superalgebra on a finite set of even and odd generators. 

\vskip .2cm
 (d) The de Rham algebra $\Omega^\bullet_X$ of differential forms on a $C^\infty$-manifold $X$, is a $\ZZ$-graded supercommutative
 algebra over $\RR$.
 
 \vskip .2cm
 
 (e) The cohomology algebra $H^\bullet(X,\k)$ of any CW-complex $X$, is a $\ZZ$-graded supercommutative
 algebra over $\k$. 

\end{exas}

The importance and appeal of such studies are based on the following heuristic

\begin{para}\label{para:nat-super}
{\bf Principle of Naturality   of Supers.}
   All constructions and features that make commutative algebras special among all algebras, can be extended
   to
supercommutative algebras, and make them just as special.  

\end{para}

The most important of such features is the relation to geometry: a commutative algebra $A$ can be seen as the algebra
of functions on a geometric object $\Spec(A)$ which can be constructed from $A$ and used to build more complicated
geometric objects by gluing. Supergeometry (as understood by mathematicians)
is the study of similar geometric objects for supercommutative algebras.


\subsection {The symmetric monoidal category of super-vector spaces}

The main guiding principle for extending properties of usual commutative algebras to supercommutative ones is also
sometimes called the Koszul sign rule. It goes like this:

\begin{para}\label{para:koszul}
  When we move any quantity (vector, tensor, operation) of parity $p\in\ZZ/2$ past any other quantity of parity $q$,
this move should be accompanied by multiplication with $(-1)^{pq}$.  
\end{para}

An instance is given by the first formula in Example \ref{exas:supercomm}(c). This rule can be formalized as follows.

\vskip .2cm

Let $\SVect_\k$ be the category of {\em super-vector spaces} over $\k$, i.e.,  $\ZZ/2$-graded vector spaces
$V=V^\ev\oplus V^\od$. This category has a  monoidal structure $\otimes$, the usual graded tensor product. The operation
$\otimes$ is associative up to natural isomorphisms, and has a unit object $\1 =\k$ (put in degree $\ev$). Define the
symmetry isomorphisms
\be\label{eq:symmetry}
R_{V,W}: V\otimes W \lra W\otimes V, \quad v\otimes w \mapsto (-1)^{\deg(v)\cdot\deg(w)} w\otimes v.  
\ee

\begin{prop}
The family $R=(R_{V,W})$ makes $(\SVect_\k, \otimes, \1, R)$ into a symmetric monoidal category. 
\end{prop}

For background on symmetric monoidal categories we refer to \cite{deligne-tannakiennes, maclane}. 
A basic example of a symmetric monoidal category is given by the usual category of vector spaces $\Vect_\k$, with
the usual tensor product, the unit object $\1=\k$, and the symmetry given by $v\otimes w\mapsto w\otimes v$. 
The meaning of the proposition is that $\SVect_\k$ with  symmetry \eqref{eq:symmetry} satisfies all the same formal
properties as the ``familiar" category $\Vect_\k$.

\vskip .2cm

If $\dim(V^\ev)=m$ and $\dim(V^\od)=n$, we write $\dim(V)=(m|n)$. In particular, we have the standard coordinate superspaces
$\k^{m|n}$. We denote by $\Pi$ the {\em parity change functor} on $\SVect_\k$ given by multiplication with $\k^{0|1}$ on the left. 

\vskip .2cm

It is well known  that one can develop linear algebra formalism (tensor, symmetric, exterior powers etc.) in any symmetric
monoidal $\k$-linear abelian category $(\Vc, \otimes, \1, R)$. 
This gives a way to define commutativity. That is, 
  a {\em  commutative algebra in $\Vc$} is an object $A\in\Vc$
together with a morphism $\mu_A: A\otimes A\to A$ satisfying associativity and such that the composition
\[
A\otimes A \buildrel R_{A,A}\over\lra A\otimes A\buildrel \mu_A\over\lra A
\]
is equal to $\mu_A$. Given two commutative algebras $(A, \mu_A)$ and $(B, \mu_B)$, the object $A\otimes B$ is again
a commutative algebra with respect to $\mu_{A\otimes B}$ given by the composition
\[
A\otimes B\otimes A\otimes B \buildrel \Id\otimes R_{B,A}\otimes \Id \over\lra A\otimes A\otimes B\otimes B
\buildrel \mu_A\otimes\mu_B \over\lra A\otimes B. 
\]
For $\Vc=\SVect_\k$, this gives Example  \ref{exas:supercomm}(c). 

\begin{rems}

(a)  Principle \ref{para:nat-super} can be expressed more formally by saying that most constructions involving commutative algebras
can be expressed in terms of the symmetry isomorphism of the tensor product and so can be reproduced in any symmetric monoidal category.
The category $\SVect_\k$ is  of course a prime example. More generally, super-generalizations of Tannakian categories studied by Deligne \cite{deligne-tensor}
provide a wider framework for extending  the formalism of commutative algebra. 

\vskip .2cm

(b) Symmetric monoidal categories can be seen as categorical analogs of commutative algebras: instead of
the equality $ab=ba$, we now have canonical isomorphisms $V\otimes W \simeq W\otimes V$. 
In Section \ref{sec:higher-super} below  we will define categorical analog of {\em super}commutative algebras. 

\end{rems}


\subsection{Superschemes and supermanifolds}

Given a supercommutative algebra $A$, the even part $A^\ev$ is commutative.
At the same any element $a$ of  the odd part $A^\od$ is nilpotent: $a^2=0$.

The commutative algebra $A^\ev$  has the associated affine scheme
$\Spec(A^\ev)$. Explicitly, 
 it is a ringed space $\bigl(\ul{\Spec}(A^\ev), \Oc_{\Spec(A^\ev)}\bigr)$, where $\ul{\Spec}(A^\ev)$ is the set of prime
ideals in $A^\ev$ with the Zariski topology and $\Oc_{\Spec(A^\ev)}$ is a sheaf of local rings on this space obtained by localization of $A^\ev$.
That is, the value of $\Oc_{\Spec(A^\ev)}$ on a ``principal open set'' 
\[
U_f = \{ \pen\in\ul\Spec(A^\ev): f\in\pen\}, \quad f\in A^\ev,
\]
is given by $\Oc_{\Spec(A^\ev)}(U_f) = A^\ev[f^{-1}]$. 

\vskip .2cm

Further, $A^\ev$ is not only commutative but lies in the center of $A$ as an associative algebra. Therefore, associating to 
$U_f$ the commutative superalgebra $A[f^{-1}]$, we get a sheaf of commutative superalgebras
\[
\Oc_{\Spec(A)} \,\,=\,\, \Oc_{\Spec(A)}^\ev \oplus \Oc_{\Spec(A)}^\od, \quad \Oc_{\Spec(A)}^\ev = \Oc_{\Spec(A^\ev)}
\]
on the same topological space $\ul\Spec(A^\ev)$ as before. The pair (ringed space) 
\[
\Spec(A) \,\, := \,\, \bigl( \ul\Spec(A^\ev), \Oc_{\Spec(A)}\bigr)
\]
is the fundamental geometric object associated to $A$, see \cite{leites, M}. 

The stalks of  $\Oc_X$ are commutative superalgebras which,
considered as ordinary associative algebras, are  {\em local rings}. Indeed, we have:

\begin{prop}
If $B$ is a commutative superalgebra such that $B^\ev$ is a local ring with maximal ldeal $\men^\ev$,
 then $B$ itself is a local ring with maximal ideal $\men = \men^\ev \oplus B^\od$.  \qed
\end{prop}

\noindent {\sl Proof:}  Odd elements of a supercommutative algebra are nilpotent. 
This implies at once that $b=b_\ev+b_\od\in B$ is invertible,  if and only if $b_\ev$ is invertible in $B^\ev$, i.e., $b\notin\men$.
In other words, $B$ is a ring in which non-invertible elements form an ideal, i.e., a local ring. \qed

\vskip .2cm

Nilpotency of odd elements mentioned above, 
 is the reason why the underlying space of $\Spec(A)$ depends only on the even part of $A$. 
 Similarly to the usual commutative algebra intuition about spectra of rings with nilpotents,
 $\Spec(A)$ may be thought of, roughly, as some kind of ``infinitesimal neighborhood thickening'' of $\Spec(A^\ev)$.

\vskip .2cm

Working  with sheaves of local rings is a fundamental technical feature of Grothendieck's theory of schemes. So we introduce 
the following

\begin{Defi} 
A   {\em super-locally ringed space over $\k$}  is a pair 
 $X=(\ul X, \Oc_X)$, where $X$ is a topological space and $\Oc_X=\Oc_X^\ev\oplus \Oc_X^\od$
is a sheaf of commutative $\k$-superalgebras on $\ul X$, with stalks being local rings. 

A {\em morphism} of super-locally ringed spaces
$f: (\ul X,\Oc_X)\to (\ul Y, \Oc_Y)$ is a pair $f=(f_\sharp, f^\flat)$, where:

\begin{itemize}  
\item $f_\sharp: \ul X\to \ul Y$
is a continuous map of topological spaces.

\item  $f^\flat: f_\sharp^{-1}(\Oc_Y)\to\Oc_X$ is a morphism of sheaves
of commutative superalgebras, which, in addition, takes the maximal ideal of each local ring
$(f_\sharp^{-1}\Oc_Y)_x : = \Oc_{Y, f_\sharp(x)}$,  $x\in X$ into the maximal ideal of the local ring $\Oc_{X,x}$. 
\end{itemize}
The resulting category of super-locally ringed spaces  over $\k$ will be denoted $\SLRS$

\end{Defi} 

\begin{prop}
For any two commutative superalgebras $A,B$ we have an isomorphism
\[
\Hom_{\SA}(A,B) \,\,\simeq \,\,\Hom_{\SLRS}(\Spec(B), \Spec(A)). 
\]

\end{prop}

\begin{Defi}
A {\em superscheme} over $\k$ is a super-locally ringed space over $\k$ locally isomorphic to $\Spec(A)$ for a commutative
superalgebra $A$. The category $\SSch$ of superschemes over $\k$ is defined to be the full subcategory of
$\SLRS$ formed by superschemes. 
\end{Defi}

\begin{exas}  (a)  The {\em affine superspace} of dimension $(m|n)$ is the superscheme
\[
\AAA^{m|n}\,\,=\,\,\Spec \bigl( \k[x_1, \cdots, x_m]\otimes\Lambda[\xi_1, \cdots, \xi_n]\bigr).
\]

(b) 
The next class of examples is provided by {\em algebraic supermanifolds} of dimension $(m|n)$. By definition, they are 
superschemes locally of the form $\Spec \bigl( A\otimes \Lambda[\xi_1, \cdots, \xi_n]\bigr)$, where $A$ is the ring of
functions on a smooth $m$-dimensional affine  algebraic variety over $\k$. 
\end{exas}

Each superscheme $X=(\ul X, \Oc_X)$ gives rise to an ordinary scheme $X^\ev=(\ul X, \Oc_X^\ev/(\Oc_X^\od)^2)$
called the {\em even part} of $X$. If $X$ is an algebraic supermanifold of dimension $(m|n)$, then
$X^\ev$ is a smooth algebraic variety of dimension $m$. 

As with usual schemes, a superscheme $X$ is defined by its {\em functor of points} on the category of supercommutative algebras
\[
A \,\,\mapsto \,\, X(A) \,=\,\Hom_{\SSch} (\Spec(A), X). 
\]

Differentiable or analytic supermanifolds are similarly defined as locally ringed spaces $X=(X^0, \Oc_X)$ where $X^0$
is an ordinary differentiable or analytic manifold and $\Oc_X$ is a sheaf of commutative superalgebras locally isomorphic
to $\Oc_X\otimes_\k \Lambda[\xi_1, \cdots, \xi_n]$. Here $\Oc_X$ is the sheaf of $C^\infty$ or analytic functions on $X$
and $\k=\RR$ for differentiable or real analytic and $\CC$ for complex analytic manifolds.

 
 \subsection{Lie supergroups and superalgebras}
 
 Here is an illustration of the Principle of Naturality of Supers: extension of the formalism of differential geometry
 to the super-case looks totally straightforward. 
 
 Given a $\k$-linear symmetric monoidal category $\Ac$,
 one can speak about algebras in $\Ac$  of any given type, for example Lie algebras. In the case $\Ac=\SVect$,
this gives the familiar concept.

\begin{Defi}
A {\em Lie superalgebra} over $\k$ is a $\k$-super-vector space $\gen = \gen^\ev\oplus\gen^\od$ equipped with a  homogeneous
(degree $\ev$) operation $(x,y)\mapsto [x,y]$ satisfying
\[
 \beg
 \, [x,y] = 
 - (-1)^{\deg(x)\cdot \deg(y)} [y,x]; 
\\
 [x,[y,z]] \,\,\,=\,\,\, [[x,y],z] \, + \, (-1)^{\deg(x)\cdot \deg(y)} [y,[x,z]].
\end{gathered}
\]
\end{Defi}

The following examples can also  be given in any symmetric monoidal category, we use the case $\Vc=\SVect_\k$. 
\begin{exas}
(a) If $R$ is any associative superalgebra, we can make it into a Lie superalgebra using the {\em supercommutator}
\[
[x,y] \,\,=\,\, x\cdot y - (-1)^{\deg(x)\cdot \deg(y)} y\cdot x.
\]
(b) If $A$ is an associative superalgebra, then we have the Lie superalgebra $\Der(A)=\Der^\ev(A)\oplus\Der^\od(A)$
of {\em super-derivations} of $A$. Here
 $\Der^i(A)$ is the space of $\k$-linear maps $D: A\to A$ satisfying the {\em super-Leibniz rule}
\[
D(a\cdot b) = D(a) \cdot b + (-1)^{i\cdot \deg(a)} a\cdot D(b). 
\]
The bracket is given by the supercommutator as in (a). 
\end{exas}

Given a supermanifold $X=(\ul X,\Oc_X)$  of dimension $(m|n)$ in the smooth, analytic or algebraic category, super-derivations
of $\Oc_X$ form a sheaf  $T_X$ 
of Lie superalgebras on $X$ called the {\em tangent sheaf}.
 It is also a sheaf of $\Oc_X$-modules,  locally free of rank $(m|n)$. Its fiber at a  point
$x\in X^0$ will be denoted $T_xX$ and called the {\em tangent space to $X$ at $x$}.
 It is a super-vector space of dimension $(m|n)$. 

\vskip .2cm

A {\em Lie supergroup}, resp. an {\em algebraic supergroup} over $\k$, is a group object $G$ in the category of
smooth supermanifods, resp. algebraic supermanifolds over $\k$. In particular, it has a unit element $e\in G^0$.
The tangent space $ T_e G$ is a Lie superalgebra denoted by $\gen=\Lie(G)$. Note that $\gen^\ev=\Lie(G^\ev)$
is the usual Lie algebra corresponding to a Lie or algebraic group. 
 Conversely, given a finite-dimensional Lie superalgebra $\gen$
and a Lie or algebraic group $G^\ev$ such that $\Lie(G^\ev)=\gen^\ev$, one can integrate it to a Lie  or algebraic supergroup
$G$ with $\Lie(G)=\gen$.

We will specially use the case when $\gen$ is {\em nilpotent}, that is, the iterated commutators
$[x_1, [x_2, \cdots, x_n]...]$ vanish for $n$ greater than some fixed $N$ (degree of nilpotency). 
We assume $\on{char}(\k)=0$. 
In this case we can associate  to $\gen$ a canonical  algebraic  group $G=e^\gen$
by means of the
classical Hausdorff series
 \be\label{eq:hausdorff}
 x\cdot y \,\,=\,\, x+ y + {1\over 2} [x,y] + \cdots
 \ee
  This means the following. Suppose we want to find the set-theoretic group of $A$-points
   $e^\gen(A)=\Hom(\Spec(A), e^\gen)$,
 where $A$ is a super-commutative algebra. As a set, this group is defined to be 
 $(\gen\otimes_\k A)^\ev$, which has a structure of an {\em ordinary} (purely even) nilpotent  Lie algebra and therefore
 can be made into a group by means of the {\em ordinary} Hausdorff series above. This is the group $e^\gen(A)$. 
 
 Informally, one says that $e^\gen$ ``is'' the Lie algebra $\gen$ considered as a manifold with the multiplication given by
 \eqref{eq:hausdorff} but the  formal meaning is   that the  arguments in commutators in  \eqref{eq:hausdorff} always
 belong to some ordinary Lie algebra. 
 
 \begin{ex}[(The super-particle and its automorphism group)] In classical mechanics, a ``particle'' is thought of as a point, i.e., mathematically, as the
 the variety $\Spec(\k)$. It is natural to call the {\em super-particle} the affine superspace $\AAA^{0|1} = \Spec( \Lambda[\xi])$.
 Although  has only one classical ($\k$-)point, the super-structure around this point is different.  As shown by Witten
 \cite{witten-morse}, study of such a super-particle moving in an (ordinary) Riemannian manifold leads naturally
 to    the de Rham complex $\Omega^\bullet_M$  and to Morse theory.

 In particular, we have
 the algebraic supergroup $\Aut(\AAA^{0|1})$ formed by ``all automorphisms" of $\AAA^{0|1}$, even as well as odd.
  For a commutative $\k$-superalgebra $A$ the group of $A$-points
 of  $\Aut(\AAA^{0|1})$ consists of $A$-super-algebra automorphisms of $A\otimes_\k\Lambda[\xi]$, i.e., of transformations 
 \[
 \xi\,\,\mapsto\,\,  a\xi+b, \quad a\in A_\ev^*\,  \text{ (invertible)}, \,\,\, b\in A_\od. 
 \]
 This means that $\Aut(\AAA^{0|1})$ has dimension $(1|1)$ and its even part is the multiplicative group $\GG_m$. 
 The following proposition is due to M. Kontsevich \cite{kontsevich}, see \cite{KV} Prop. 2.2.2 for a discussion.
 It explains the appearance of the de Rham differential in Witten's theory. 
 
 \begin{prop}\label{prop:super-dg}
 Let $V$ be a super-vector space over $\k$. An (algebraic) action of $\Aut(\AAA^{0|1})$ on $V$ is the same as the following
 system of data:
 \begin{itemize}
 \item[(1)] A $\ZZ$-grading $V=\bigoplus_{in\in \ZZ} V^i$ such that $V_\ev = \bigoplus V^{2i+1}$ and $V_\od = \bigoplus V^{2i}$. 
 
 \item[(2)] A differential $d$ on $V$ of degree $+1$ in the above $\ZZ$-grading such that $d^2=0$, i.e., making $V$ into
 a cochain complex. 
 \end{itemize}
 In other words, the symmetric monoidal category of   super vector spaces with an $\Aut(\AAA^{0|1})$-action is  identified
  with the symmetric monoidal category
 of  cochain complexes.  \qed
 \end{prop}
 
 Here the $\ZZ$-grading is given by the action of $\GG_m\subset \Aut(\AAA^{0|1})$, and the differential is given by
 the action of the odd part of the Lie super-algebra $\Lie(\Aut(\AAA^{0|1}))$ of super-derivations of $\Lambda[\xi]$. 
  \end{ex}
  
  \begin{parab}{\bf Super-geometry and derived geometry.} It is instructive to compare super-geometry (as it is described in
  this section) with the formalism of {\em derived geometry} which seeks to form ``nonabelian derived functors''
  of several familiar algebro-geometric constructions such as forming moduli spaces or intersections of subvarieties,
  see \cite{CFK} \cite{lurie} \cite{TV}.

  In the simplest setting,  objects studied  in derived geometry are glued out of local pieces which correspond to
  commutative differential graded (dg-) algebras $(A^\bullet, d)$, i.e., commutative algebra objects in the symmetric monoidal
  category of cochain complexes.  By Proposition \ref{prop:super-dg},
  such pieces can be seen as affine super-schemes with action of   $\Aut(\AAA^{0|1})$.
  
  So one can say that, in some approximate sense, derived geometry can be  regarded as super-geometry but in an
  equivariant context, with respect to the action of the supergroup  $\Aut(\AAA^{0|1})$. However, this is a rather
  simplified point of view for several reasons. Even in the context of dg-algebras, the distinction between 
  negative and positive grading is very important, it corresponds to the distinction between ``geometry in the small''
   (intersections, singularities)  and  ``geometry in the large'' (stacks, homotopy type). The precise definitions impose 
   ``geometry in the large''
  in a  more global ``stacky'' way.

  \end{parab}
 
 \begin{parab}\label{parab:pfaff}
 {\bf Pfaff systems and Frobenius theorem.} Let $X$ be a supermanifold. A {\em Pfaff system} in $X$ is a subbundle (i.e., a subsheaf
 of $\Oc_X$-modules 
 which is locally a direct summand) $C\subset T_X$. For a Pfaff system $C$ the Lie algebra structure in sections of $T_X$
 induces the $\Oc_X$-linear map known as the {\em Frobenius pairing}
 \[
\digamma_C:  \Lambda^2_{\Oc_X} C \lra T_X/C.
 \]
 A Pfaff system $C$ is called {\em integrable}, if $\digamma_C=0$, i.e., $C$ is closed under the bracket of vector fields.
 In this case we a super-analog of the {\em Frobenius theorem}: in the smooth or analytic case $C$ can be locally represented
 as the relative tangent bundle to a submersion of supermanifolds $X\to Y$. 
 \end{parab}


 \section{Supergeometry as understood by physicists}
 
 For a physicist, the really important concept is {\em supersymmetry}, and supermanifolds per se are of interest only
 tangentially. For example, they arise as infinite-dimensional spaces of bosonic and fermionic classical fields,
 over which Feynman integrals are taken. One way (not the only one!) to construct supersymmetric field theories is
 by using {\em superspace}, a concept not synonymous with ``supermanifold'' of mathematicians.  In fact, superspace
 is a supermanifold with a rather special ``spinor-conformal" structure.

 \subsection {The idea of non-observable square roots} 
 To understand the idea of supersymmetry, it is useful to include it into the following more general heuristic principle.
 
 \begin{para}{\bf Principle of square roots.} It is useful to represent  observable quantities of immediate physical interest
 (real, positive, bosonic) as bilinear combinations of more fundamental quantities which can be complex, fermionic and
 not even  observable by themselves. 
 \end{para}
 
In other words, it is useful to take ``square roots'' of familiar objects. Let us give several examples.

 \begin{ex}[(Wave functions and probability density)]\label{ex:wave}
 In elementary quantum mechanics, the wave function $\psi(x)$ of a particle
(say, electron), is a complex quantity which can not be measured. But the expression
\[
P(x) \,\,=\,\, |\psi(x)|^2 \,\,=\,\, \ol\psi(x) \cdot \psi(x) \,\, \geq 0
\]
represents the probability density of the electron which is real, non-negative and measurable.

 \end{ex}
 
 \begin{ex} [(Laplace operator on forms)]\label{ex:laplace}
 Let $X$ be a $C^\infty$ Riemannian manifold. The space $\Omega^\bullet(X)$ of all differential forms on $X$
 is $\ZZ$-graded and so can be considered as a super-vector space. The Laplace operator on forms
 is defined as
 \[
 \Delta \,\,=\,\, d\circ d^* + d^* \circ d \,\,=\,\, [d,d^*],
 \]
 where $d$ is the exterior derivative and $d^*$ is its adjoint with respect to the Riemannian metric. Thus $\Delta$
 is non-negative definite, real (self-adjoint) and bosonic, while $d$ and $d^*$ are fermionic. 
 \end{ex}
 
 \begin{ex}
 [(Spinors as square roots of vectors)]
 Let $V$ be a $d$-dimensional vector space over $\k = \RR$ or $\CC$, equipped with a non-degenerate
 quadratic form $q$. We refer to $V$ as the spacetime, or the Minkowski space.
 Elements of $V$ ({\em vectors} in the physical sense of the word)
  can be represented as bilinear combinations of more fundamental
 quantities, {\em spinors}. By definitions, spinors are vectors in minimal (spinor) representations of $\on{Spin}(V)$, 
 a double covering of the group $SO(V)$. There is one spinor representation $S$, if $d=\dim(V)$ is odd, and two representations,
 $S_+$ and $S_-$, if $d$ is even. Further properties of these representations depends on the residues modulo 8
 of $d$ and of the signature of $q$, if $k=\RR$ (Bott periodicity). We refer to \cite{ABS, deligne-spinors} for a detailed treatment. 
 
 The expression of vectors as bilinear combinations of spinors comes from $\Spin(V)$-equivariant maps $\gamma$
 whose nature we indicate in three cases $d=2,4,10$ with Minkowski signature.
 
 \begin{itemize}
 \item[$d=2$:]  $S_+, S_-$ can be defined over $\RR$ and have real dimension 2. They are self-dual: $S_\pm^*=S_\pm$. The
 gamma maps are $\gamma: \Sym^2(S_\pm)\to V$.
 
 \item[$d=4$:] $S_+, S_-$ are complex and have complex dimension 2. They are hermitian conjugate of each other: $\ol S_\pm^*=S_\mp$.
 The gamma map is $\gamma: S_+\otimes S_- \to V$.
 
 \item[$d=10$:] $S_+, S_-$ are real, of dimension 16 and self-dual. The gamma maps are $\gamma:\Sym^2(S_\pm)\to V$. 
 
 \end{itemize}
  \vskip .2cm
 
 The maps $\gamma$   can be seen as transposed versions of   the systems of Dirac's gamma-matrices. For instance, in the case $d=10$,
we can transpose $\gamma$ to get a map $\gamma^t: V\to \Hom(S_\pm,  S_\pm)$, and the gamma-matrix $\gamma_\mu$,
$\mu=1,\cdots,  10$,  is the operator
$\gamma^t(e_\mu): S_\pm\to S_\pm$ corresponding to the basis vector $e_\mu\in V$.

 \vskip .2cm
 
  Needless to say,  forming the double covering $\Spin(V)\to SO(V)$ itself can be seen as taking square roots of rotations. The appearance
 of spinors is a reflection of this procedure at the level of representations of groups. 
 
 \end{ex}

 \begin{ex}
 [(Weil conjecture over finite fields)] It is tempting to add to this list of ``square roots'' the following classical example. Let $X$ be
 a smooth projective curve over a finite field $\FF_q$. The \'etale cohomology group $H^1(X\otimes\ol\FF_q, \QQ_l)$
 is acted upon by the Frobenius element $\on{Fr}$, generating $\on{Gal}(\ol\FF_q/\FF_q)$. It was proved by A. Weil
 (and motivated the more general Weil conjectures) that each eigenvalue $\lambda$ of $\on{Fr}$ is an algebraic
 integer whose image in each complex embedding satisfies $|\lambda| = \sqrt{q}$, i.e., $\lambda\cdot\ol\lambda = q$.
 So the first cohomology, a fermionic structure, gives rise to factorizations $q= \lambda\cdot\ol\lambda$. 
 
 Further, it was suggested by Y. I.  Manin that the motive of a supersingular elliptic curve over $\FF_q$ can be
 seen as a ``spinorial square root" of the Tate motive. See \cite{rama} for developments in this direction. 
 \end{ex}
 
 \begin{rem}
 Example \ref{ex:wave} (understood more widely,  as an instance of the special role played  by scalar products and Hilbert
 spaces in quantum mechanics),  can be seen as 
 the  source of many appearances of square roots.
  For instance, a natural $L_2$-space associated to a manifold $M$,  consists of half-forms, i.e., of
 sections of  $\omega_M^{\otimes (1/2)}$, a { square root} of the line bundle of volume forms. 
 
  The appearance of the metaplectic double cover $\wt{\on {Sp}}(n)\to \on{Sp}(n)$   of the symplectic group can also be traced, via the Maslov index and the WKB
  approximation, to the same source.

 \end{rem} 
 
 \begin{ex}[(Connections and curvature)]\label{ex:connections}
 A gauge (say, electromagnetic) field  is represented by a connection $\nabla$ in a principal bundle
 over the space-time manifold. The  immediately observable quantity (at least  for the electromagnetic field which  can be measured classically)
 is  the {\em field strength}, represented by the curvature  $F=\nabla^2$. 
 The connection itself, a somewhat more elusive   object,  can be viewed therefore as a square root of the field strength. 
 The idea that ``the fundamental quantities that we observe, are curvature data of something''  appears, of course, in several
 areas  of physics.  Its subtle match with the basic principles of quantum mechanics is a truly remarkable phenomenon.

 \end{ex}
 
 
 \subsection{ A square root of $d/dt$ and theta-functions}\label{subsec:theta}
 The simplest example of supersymmetry can be obtained by looking  at the differential operator (super-derivation)
 \be\label{eq:root}
 Q \,\,=\,\, {\partial\over\partial \xi} + \xi {\partial\over\partial t} \,\,=\,\, \begin{pmatrix}
 0&1
 \\
 \partial/\partial t & 0
 \end{pmatrix}, \quad Q^2={\partial\over\partial t}
 \ee
 in the algebra $\Oc(\AAA^{1|1}) = \CC[t] \otimes\Lambda[\xi]$. One can replace $\CC[t]$ by the algebras of smooth or analytic
 functions of $t$. The second equality sign in \eqref{eq:root} is an instance of {\em component analysis}: 
 we view $\CC[t]\otimes\Lambda[\xi]$ as a free $\CC[t]$-module with basis $1,\xi$ and write $Q$ as a $2\times 2$ matrix
 differential operator in $t$ alone. One verifies immediately that $Q^2={\partial/ \partial t}$, so $Q$ gives a square root
 of the Hamiltonian (time translation). 
 
 Already this example is quite non-trivial: it  provides a natural explanation of Sato's approach to theta functions 
 \cite{Sato-theta, SKK-theta}.
 We do it in two stages. First, consider the exponential
 \be
 e^Q \,\,=\,\, 1+ Q + {Q^2\over 2!} + {Q^3\over 3!} + \cdots
 \ee
 as series of operators acting in complex analytic functions in $t,\xi$. 
 
 \begin{prop}\label{prop:eQ}
 $e^Q$ converges to   a local operator, i.e., to an endomorphism of the sheaf $\Oc_{\CC^{1|1}}$
 of analytic functions on $\CC^{1|1}$.
 \end{prop}
 
 A standard example of a local operator in this sense is 
  a  linear differential operator $P$  on the line.   An   operator with constant coefficients
   is just a polynomial in $d/dt$:
\[
P \,\,=\,\, h(d/dt), \quad h(z) = \sum_{n=0}^N a_n z^n.
\]
Replacing polynomials by entire analytic functions $h(z)=\sum_{n=0}^\infty a_n z^n$, we get expressions which, even when they converge, need not give local
operators. For instance, $h(z)=e^z$ gives the shift operator
\be\label{eq:shift}
(e^{d/dt} f)(t) \,\,=\,\, f(t+1). 
\ee
 However, if $h(z)$ is {\em sub-exponential}, i.e., the series $k(z) = \sum n! a_n z^n$ still represents an entire function,
 then $h(d/dt)$ acts on analytic functions in a local way. Indeed, by  the Cauchy formula,
 \[
 \bigl(h(d/dt)f\bigr)(t) \,\,=\,\,\sum_{n=0}^\infty a_n f^{(n)}(t) \,\,=\,\,{-1\over 2\pi i} \oint_{|t'-t|=\eps}
 f(t') k\biggl({1\over{t-t'}}\biggr) dt',
 \]
 where $\eps$ can be arbitrarily small. So if $f$ is analytic in an open $U\subset \CC$, then so is $h(d/dt)f$.
 For example, $h(z) = \cos \sqrt{z}$ gives a local operator. One can similarly make sense of series $\sum_{n=0}^\infty a_n(t) (d/dt)^n$,
 where $a_n(t)$ are analytic functions in $t$ of  sub-exponential growth in $n$. These series define local operators on
 analytic functions on $\CC$ known simply as {\em differential operators of infinite order} \cite{SKK-theta}. They form a sheaf
 of rings on $\CC$, denoted $\Dc_\CC^\infty$. Similarly for any complex analytic (super-)manifold such as $\CC^{1|1}$. 
 
  \vskip .2cm
  
  Returning to the situation of Proposition \ref{prop:eQ},  we see that $Q$, being a square root of $\partial/\partial t$, has exponential
  $e^Q$ which is  a local operator:  a global section of $\Dc_{\CC^{1|1}}^\infty$ or, after component analysis,
  of $\on{Mat}_2(\Dc_\CC^\infty)$. 
  
  \begin{rem}
  Although $Q$ is a super-derivation, it is not a derivation in the usual sense and  $e^Q$ is not a ring automorphism.
  In particular, forming $e^Q$ is not an instance of exponentiating a Lie superalgebra to a Lie supergroup.

  \end{rem}
  
  \vskip .2cm
  
  We now consider the simplest theta-function
  \[
  \theta(t,x) \,\,=\,\,\sum_{n\in\ZZ} e^{n^2t + nx}, \quad \Re(t)<0, \, x\in\CC. 
  \]
 Its value at $x=0$ (the Thetanullwert)
 \[
 \theta(x,0) \,\,=\,\,\sum_{n\in\ZZ} q^{n^2}, \quad q=e^t, \,\, |q|<1,
 \]
 is a modular form, so there is a relation (modular transformation) relating its values at $t$ and at $1/t$ (in our normalization). 
 The main result of \cite{SKK-theta} is a characterization of $\theta(t,0)$ by two differential equations of infinite order in $t$ alone
 (i.e., by local conditions in $t$).   They are then used to deduce
  the modular transformation  because the system of  
 equations is invariant under it. 
 
 This can be done as follows. Note that $\theta(t,x)$, as a function of two variables, satisfies the heat equation
 \[
 {\partial\theta\over\partial t} \,\,\,=\,\,\, {\partial^2\theta\over\partial x^2}. 
 \]
 This means that $\partial/\partial x$ acts on $\theta$ as a square root of $\partial/\partial t$. On the other hand, $\theta$
 has periodicity properties
 \[
 \begin{cases}
 \theta(t, x+2\pi i) \,\,=\,\,\theta(t,x)
 \\
 \theta(t, x+2t) \,\,=\,\, e^{-x-t} \theta(t,x).
 \end{cases}
 \]
Using \eqref{eq:shift} in the $x$-variable and  the heat equation, we can write this formally as
 \[
 \begin{cases}
 e^{2\pi i \sqrt{\partial/\partial t}} \theta(t,x) \,\,=\,\,\theta(t,x)
 \\
 e^{2 t \sqrt{\partial/\partial t}}\theta(t,x) \,\,=\,\, e^{-t-x} \theta(t,x),
 
 \end{cases}
 \]
 after which we can specialize to  $x=0$. 
 Now, replacing $\sqrt{\partial/\partial t}$ with $Q$, we get a system of two differential equations of infinite order in $t$ alone,
 satisfied by $\theta(t,0)$.

 
 \subsection{Square roots of spacetime translations}
 Representing just the Hamiltonian (the operator of energy, or time translation)
 as a bilinear combination of fermionic operators (\S \ref{subsec:theta}, or Example \ref {ex:laplace}), is a non-relativistic procedure.
  Relativistically, 
 we cannot separate $\partial/\partial t$ from any other constant vector field $\partial_v$ (momentum operator) on the
 Minkowski space $V$ and so should represent them all. This leads to the following concept. 
 
 \begin{Defi}
 A {\em quadratic space} over a field $\k$ is a datum of finite-dimensional $\k$-vector spaces $V,B$ and a surjective
 linear map $\Gamma: \Sym^2(B)\to V$. We will often view  $\Gamma$ as a $V$-valued scalar product
 $(b_1, b_2)\mapsto \Gamma(b_1, b_2)\in V$ on $B$
  \end{Defi}
 
 A quadratic space $\Gamma$ gives rise to the Lie superalgebra
 \be
 \ten = \ten_\Gamma,\quad \ten^\ev= V, \,\,\, \ten^\od = B,
 \ee
 with the only non-zero component of the bracket being $\Gamma$.  We call $\ten$ the {\em supersymmetry algebra}
 associated to $\Gamma$. The abelian central subalgebra $\ten^\ev=V$ is the usual Lie algebra of infinitesimal spacetime
 translations, and we denote by $T^\ev=e^{\ten^\ev}$ the corresponding algebraic group (i.e., $V$ considered as an algebraic
 variety and equipped  with vector addition). Since $\ten$ is nilpotent: $[x,[y,z]]=0$ for any homogeneous $x,y,z$, it is easily integrated
 to an algebraic supergroup $T=e^\ten$ called the {\em supersymmetry group} using the truncated Hausdorff series
 \eqref{eq:hausdorff}.

 \begin{exas}\label{ex:susy}
 In physical applications (see \cite{deligne-freed, Freed-super} for a detailed exposition), $V$ is the Minkowski space 
 (i.e., a vector space over
 $\RR$ or $\CC$ with a  non-degenerate quadratic form) and $B$ is a direct sum of (possible several copies of)
 the spinor bundle(s). Different possible choices of such $B$ are known as different types of 
 {\em extended supersymmetry}. More precisely:
 
 \begin{itemize}
 \item[] $N=p$ {\em supersymmetry (SUSY)} means that $d=\dim(V)$ is odd and $B=S^{\oplus p}$.
 
 \item[] $N=(p,q)$ {\em supersymmetry (SUSY)} means that $d$ is even and $B= S_+^{\oplus p} \oplus S_-^{\oplus q}$. 
 
 \end{itemize}
 The map $\Gamma$ of the quadratic space is constructed out of the gamma-maps $\gamma$ for spinors.
 Quadratic spaces obtained in this way will be called {\em spinorial}. 
 For example:
 
 \vskip .2cm
 
 (a) $d=10$, $N=(1,0)$ SUSY: Here $B=S_+$ and $\Gamma=\gamma: \Sym^2(S_+)\to V$.
 This is the most fundamental example in many respects. 
 
 \vskip .2cm
 
 (b) $d=4$, $N=(1,1)$ SUSY: Here $B=S_+\oplus S_-$ and $\Gamma$ is the composition
 \[
 \Sym^2(S_+\oplus S_-) \lra S_+\otimes S_- \buildrel\gamma\over \lra V. 
 \]
 
 (c) $d=2$, $N=(1,1)$ SUSY is often written explicitly using two operators $H$ (energy $=\partial/\partial t$) and
 $P$ (momentum $=\partial/\partial x$). The space $S_+$ is spanned by two vectors $Q_+, Q_+^*$ and $S_-$
 by $Q_-, Q_-^*$, and the supersymmetric algebra (i.e., the commutation rule in $\ten$) is written as
 \[
 [Q_+, Q_+] = H+P, \quad [Q_-, Q_-] = H-P, \quad [Q_+, Q_-] = 0. 
 \]
 
 For a spinorial quadratic space $\Gamma$ the group $\Spin(V)$ acts on $\ten$ and on $T=e^\ten$.
 The corresponding semidirect product $\Pc= \Spin(V) \ltimes T$ is known as the {\em super-Poincar\'e group}
 (corresponding to the type of extended supersymmetry represented by $\Gamma$).

 \end{exas}

 
 \subsection{Quadratic spaces and intersections of quadrics}
 
 In algebraic geometry, the term {\em intersection of quadrics} 
 \cite{reid, tyurin} means one of two closely related objects:
 
 \begin{itemize}
 
 \item[(1)] A {\em homogeneous intersection of quadrics}:  a subscheme $Y$ in a vector space $B$, given by
 homogeneous quadratic equations. That is, the ideal $I(Y)\subset \Sym^\bullet(B^*)$ is generated by
 $I_2(Y)\subset \Sym^2(B^*)$, its degree 2 homogeneous part.
 
 \item[(2)] A {\em projective intersection of quadrics}:  the projectivization $Z =\PP(Y)\subset  \PP(B)$
 of $Y$ as above.  In this way, projective intersections of quadrics   $Z\subset \PP(B)$ are in bijection with those homogeneous
 intersections of quadrics  $Y\subset B$ which are not entirely supported at 0.   
 \end{itemize}

 Each quadratic space $\Gamma:\Sym^2(B)\to V$ gives a homogeneous intersection of quadrics
 \[
 Y_\Gamma \,\,=\,\, \{ b \, | \, \,\Gamma(b,b)=0\} \,\,\subset \,\, B,
 \]
 and we denote by $Z_\Gamma\subset \PP(B)$ its projectivization. Alternatively,
 \[
 Y_\Gamma \,\,=\,\, \{ b\in \ten_\Gamma^\od\, | \,\, [b,b]=0\}
 \]
 is the {\em scheme of Maurer-Cartan elements} of the supersymmetry algebra $t_\Gamma$. 
 
 \vskip .2cm
 
 \begin{prop}
 Each homogeneous intersection of quadrics $Z\subset B$ can be obtained in this way from some quadratic space
 $\Gamma:\Sym^2(B)\to V$, defined uniquely up to an isomorphism.
 \end{prop}
 
 \noindent {\sl Proof:} 
 Take $V=I_2(Y)^*$, the space dual to the space of quadratic equations of $Y$ and take for $\Gamma$
 the canonical projection.
 
  \qed
 
  So  quadratic spaces are simply data encoding intersections of quadrics. 
  Classically, the simplest intersections of quadrics are as follows. 

\begin{Defi}
 (1) A projective  intersection of quadrics   $Z\subset B$ 
  is called a {\em complete
 intersection}, if $\dim I_2(Z)$, the number of quadratic equations of $Z$, is equal to the codimension of $Z$ in $\PP(B)$.
 
  \vskip .1cm
  
  (2)  A quadratic space $\Gamma:\Sym^2(B)\to V$ is called of {\em complete intersection type}, if $Z_\Gamma$
 is a complete intersection, i.e., $\dim(V)=\on{codim} (Z_\Gamma)$. 

\end{Defi}

 \begin{exas}\label{ex:susyint}
 (a) For the spinorial quadratic space $\Gamma: \Sym^2(S_+)\to V$ of $d=10$, $n=(1,0)$ SUSY (Example \ref{ex:susy}(a)),
 $X_\Gamma$ is the 10-dimensional {\em space of pure spinors} $\Sigma^{10}\subset \PP(S_+) = \PP^{15}$.
  It can be identified with
 one component of the variety of isotropic 5-planes in $V$. It is {\em not} a complete intersection: it
 has codimension 5 but is given by $d=10$
 equations.

 \vskip .2cm
 
 (b) For the quadratic space $\Gamma: \Sym^2(S_+\oplus S_-)\to V$ of $d=4$, $N=(1,1)$ SUSY (Example \ref{ex:susy}(b)),
 $X_\Gamma\subset \PP(S_+\oplus S_-) = \PP^3$ is the disjoint union of two skew lines
 $\PP(S_+) \sqcup \PP(S_-) \simeq \PP^1\sqcup \PP^1$.  It is {\em not} a complete intersection: it
 has codimension 2 but is given by $d=4$
 equations. 
 
 \vskip .2cm
 
 (c) The variety $\Sigma^{10}\subset \PP^{15}$ is a particular case of the following: a partial flag variety $G/P$
 ($G$ reductive algebraic group, $P\subset G$ parabolic), equivariantly  embedded into $\PP(E)$, where $E$ is an
 irreducible highest
 weight representation of $G$.   All such varieties are known
 to be intersections of quadrics.

 \end{exas}
 
So our physical spacetime is really {\em the space of equations} of an auxiliary intersection
(typically, not a complete intersection!) of quadrics. 

\vskip .2cm

 We will also need families of quadratic spaces parametrized by superschemes.
 
 \begin{Defi}\label{def:quadmod}
 Let $X$ be a superscheme. A {\em quadratic module} over $X$ is a datum of two locally free sheaves $B, V$ of $\Oc_X$-modules,
 both of purely even rank, and of  an $\Oc_X$-linear map $\Gamma: \Sym^2_{\Oc_X}(B)\to V$  with two properties:
 \begin{enumerate}
 \item[(1)] $\Gamma$ is surjective, and therefore the dual map $\Gamma^\vee: V^\vee\to\Sym^2_{\Oc_X}(B^\vee)$ is an embedding
 of a locally direct summand. 
 
 \item[(2)] The $\Oc_X$-algebra $\Sym^\bullet_{\Oc_X}(B^\vee)/(\Gamma^\vee(V^\vee))$ is flat over $\Oc_X$. 
 \end{enumerate}
 \end{Defi}
 
 In other words, a quadratic module gives a family of intersections of quadrics, parametrized by $X$, and
 the condition (2) means that this family is flat, in particular, its fibers have the same dimension. 
 This is important for non-complete intersections.

 
 \subsection{Supersymmetry, superspace and constraints}
 
 We start with a quick explanation of some physical terms. 
 
 \begin{parab}
 {\bf Supersymmetry} is a feature of a field theory (say, a collection of fields $\varphi$ plus a Lagrangian action $S[\varphi]$)
 defined, a priori, on the usual (non-super!) Minkowski space $V$. 
 It means that the action of the usual Poincar\'e group $SO(V)\ltimes V$ on fields by changes of variables (which leaves 
 any relativistic Lagrangian invariant)
 is extended, {\em in some way}, to an action of the super-Poincar\'e group $\Pc$ so that $S[\varphi]$ is still invariant. 
 Here $\Pc$ is constructed out of one of the spinorial quadratic spaces $\Gamma:\Sym^2(B)\to V$ (Example 
 \ref{ex:susy}). 
 
 Thus the new datum in supersymmetry is the extension of the action of $V$ to an action of $\ten_\Gamma$.
 This means that we need to 
   represent
 all the momentum operators $\partial_v, v\in V$, as bilinear combinations of  fermionic ``supercharges" $D_b, b\in B$ so that
 we have the commutation relations 
 \[
[D_b, D_{b'}] =\partial_{\Gamma(b, b')}, \quad [D_b, \partial_v]= [\partial_v, \partial_{v''}]= 0.  
 \]
 For this to be possible, there should be about equally many bosonic and fermionic fields in the theory. 
 This explains why supersymmetry is sometimes called ``symmetry between bosons and fermions''. 
  \end{parab}

\begin{parab}
{\bf  Superspace} is a tool to construct supersymmetric theories by replacing the mysterious ``in some way" above by
a natural construction. More precisely, a superspace is a supermanifold $\Sc$ extending the spacetime $V$,
(so that $V=\Sc^\ev$ is its even part) and which admits a natural action of $\ten$. 

The simplest choice
({\em flat superspace}) is $\Sc=T$, the underlying manifold of the supersymmetry group, on which $\ten = \Pi(B)\oplus V$ acts
by left-invariant vector fields $D_b, \partial_v$, see
 \cite{deligne-freed, Freed-super}.  Any field on $\Sc$ (referred to as  {\em superfield}) gives an entire multiplet of
usual fields on $V$ by {\em component analysis}: writing $\Oc_\Sc= \Oc_V \otimes\Lambda^\bullet(B^*)$
as a free $\Oc_V$-module of rank $2^{\dim(B)}$. The Lie algebra $\ten$ acts naturally on superfields, so
working only with such fields, we get supersymmetry seemingly  ``for free''

This construction can, of course, be done for an arbitrary quadratic space $\Gamma: \Sym^2(B)\to V$.
If $\Gamma$ is spinorial, then the action of $\ten$ on superfields extends to an action of the super-Poincar\'e group
$\Pc$. 

\end{parab}

 \begin{rem}
 Similarly to \S \ref{subsec:theta},   the exponentials $e^{D_b}$ of the supercharges  are local
 operators  on  analytic superfields, while the shifts $e^{\partial_v}$ are not. It would be interesting to
 understand the consequences of this phenomenon. 
 The situation of  \S \ref{subsec:theta}
 corresponds to the simplest example of a quadratic space, when $B$ is 1-dimensional,  $V=\Sym^2(B)$
 is also 1-dimensional, and $\Gamma=\Id$. 
 \end{rem}

 \begin{parab} {\bf The difficulty} with supersymmetry is that it tends to require too many fields (on $V$) for all of them
  to make physical sense. The following result  \cite{nahm} is usually intepreted by saying that ``supersymmetry in $> 11$ dimensions
 is not sensible''.
\end{parab}

\begin{para}{\bf Nahm's theorem.} Any supersymmetric theory with $d>11$ contains fields of spin $\geq 2$. 
\end{para}
 
 For the superspace construction this is easy to understand. Already the  simplest kind of  a superfield, a function on $\Sc$, is a section
 of $\Oc_V\otimes\Lambda^\bullet(B^*)$, where $B$ is the direct sum of one or several spinor spaces. As $d$ grows,
 the decomposition of $\Lambda^\bullet(B^*)$ into $\Spin(V)$-irreducibles, quickly begins to contain higher spin representations
 such as $\Sym^j(V)$, $j\geq 2$. 
 
 \vskip .2cm
 
 However, even in the remaining dimensions $d\leq 11$, the superspace construction typically gives too many component fields.
 To eliminate some of the components, one usually imposes (in a seemingly {\em ad hoc} way) some additional
 restrictions on superfields known as {\em constraints}. In the next subsection we discuss a conceptual point of
 view on such constraints.

  Supergeometry as understood by physicists, is the study of  various versions of (not necessarily flat)
  superspaces.
  All the examples that have been considered, fit into the following concept.
  
  \begin{Defi}\label{def:superspace}
  An abstract  {\em superspace} is a supermanifold $\Sc$ (smooth, analytic or algebraic)  of dimension $(m|n)$
  together with a Pfaff system  $C\subset T_\Sc$ of rank $(0|n)$ satisfying the following properties:
  \begin{enumerate}
 \item[(1)] The restriction $C|_{\Sc^\ev}$ coincides with the odd part of  $(T_\Sc)|_{\Sc^\ev}$.
  
  \item[(2)] Denote $B=\Pi(C)$ and $V=(T_\Sc)/_C$. Then, 
  the Frobenius  pairing $\digamma_C: \Lambda^2C \to (T_\Sc)/C$, written as an $\Oc_\Sc$-linear map
  $\Gamma:\Sym^2_{\Oc_{\Sc}}(B)\to V$, is a quadratic module over $\Sc$, see Definition \ref{def:quadmod}. 
   \end{enumerate}
  \end{Defi}
  
  Informally, a superspace enhancement $\Sc$ of an ordinary manifold $\Sc_\ev$ provided a framework for taking square roots of  vector fields on $\Sc_\ev$. Choice of such an enhancement can
   be viewed as a choice a square root of $\Sc_\ev$
  itself, i.e., a way of realizing the space-time  directions as some kind of curvature data, compare
  Example \ref{ex:connections}.

  \begin{ex}[(Supercurves, as understood by physicists)] For a mathematician, an (algebraic) {\em supercurve} is
  an algebraic supermanifold of dimension $1|n$ for some $n$. 
  For a physicist, a supercurve is a superspace of dimension $(1|n)$, so the Pfaff system $C$ is a necessary part of the structure.
   See \cite{DW1, DW2, M}. 
  The geometry of a supercurve of dimension $(1|1)$
  is locally modeled on the setting of
  \S \ref{subsec:theta}.  
     \end{ex}
   
   \begin{ex}[(Spinorial curved superspaces)] Definition \ref{def:superspace}  is quite general.
   It does not require that the fibers of the
   quadratic module $\Gamma$ be spinorial quadratic spaces. However, the intersections of quadrics related to spinorial spaces
   (such as the space of pure spinors) are {\em rigid} both as abstract algebraic varieties and as intersections of quadrics.
   This means that a quadratic module whose one fiber is a spinorial quadratic space, has all 
   neighboring fibers spinorial of the same type. Therefore if $\Sc$ is a superspace (in our sense)
    with the commutator pairing $\Gamma$
   spinorial at one even point, then  we have a similar  isomorphic  ``spinorial structure" in each neighboring tangent space. 
  This amounts to a differential-geometric structure on $\Sc$ including a conformal structure in the quotient bundle
  $(T_\Sc)/C$ (in particular, on the ordinary manifold $\Sc^\ev$)
   and a ``choice of spinors''  for this conformal structure.   Cf.  \cite{GGRS, lott}
   
   Cobordism categories of such curved spinorial superspaces provide a language for Atiyah-style approach to supersymmetric
   quantum field theories \cite{stolz-teichner, T}.
   
   \end{ex}
  

 \subsection{Constraints and complete intersection slices}
 
 Various recipes of imposing constrains on superfields can be  understood using the idea of 
 {\em simple plane slices of complicated intersections of quadrics}. 
 
 \vskip .2cm
 
 Let $Z\subset \PP(B)$ be an intersection of quadrics. A {\em  plane slice} of $Z$ is a scheme
 of the form
 $Z\cap M$, where   $M\subset \PP(B)$ is a projective subspace. It is an intersection of quadrics
 in $M$. The two simplest possibilities are as follows:
 
 \begin{enumerate}
 \item[(0)] $Z\cap M = Z$, i.e., $M$ is contained in $Z$ entirely.
 
 \item[(1)] $Z\cap M$ is a complete intersection of quadtrics in $M$. 
 \end{enumerate}
 
 Let $\Gamma:\Sym^2(B)\to V$ be the quadratic space corresponding to $Z$. 
 A projective subspace $M$ corresponds to a linear subspace $B'\subset B$, and $Z\cap M$
 corresponds to the   quadratic space $\Gamma': \Sym^2(B')\to V'$,
 where $V'=\Gamma(\Sym^2(B'))$ and $\Gamma'$ is the restriction of $\Gamma$. 
The supersymmetry algebra $\ten_{\Gamma'}$ is the Lie sub(super) algebra in $\ten$ generated by
 $B'\subset\ten_\Gamma^\od$.   We will call such quadratic spaces {\em slices of} $\Gamma$. 
 
 \vskip .2cm
 
 The case (0) above means that $\ten_{\Gamma'}=B'$ is abelian and purely even. 
 In the case (1) (which includes the case (0)) we will say that $\ten_{\Gamma'}$ is
 a {\em null-subalgebra}. 
 
 \begin{exas}
 (a)  Let  $Z = C\sqcup C'$ is the union of two skew lines in $\PP^3$ (Example \ref {ex:susyint}(b)).
 Clearly, each of the two lines gives an instance of Case (0) above. Other then that, 
 we   have a $\PP^1\times\PP^1$ worth of   complete intersection slices  $Z\cap M$,
 with $M$ being a chord passing through one
 point on $C$ and one point on $C'$.  
 
 \vskip .2cm
 
 (b)  Let $\Gamma: \Sym^2(S_+)\to V$ be the $d=10$ $N=(1,0)$ SUSY quadratic space. Each vector $v\in V, q(v)=0$,
 gives a linear operator $\Gamma_v: S_+\to S_+$ given by the Clifford multiplication by $v$ (transpose of the quadratic space
 map $\Gamma$). If $v$ is a non-zero null vector, i.e., $q(v)=0$, then $\Gamma_v^2=0$ and
 $\Ker(\Gamma_v)=\on{Im}(\Gamma_v)$ is an 8-dimensional subspace in $S_+$ which we denote $L_v$. 
  The intersection $\PP(L_v)\cap\Sigma^{10}$ is a quadric hypersurface in $\PP(L_v)$, and the space of equations
  of this hypersurface is spanned by $v$. In other words, $\Gamma(\Sym^2 (L_v)) = \CC\cdot v$, and
  \[
  \hen \,\,=\,\, L_v \oplus \CC\cdot  v \,\,\subset \,\, S_+\oplus V \,\,=\,\, \ten
  \]
  is a (maximal)  null subalgebra. The $(1|8)$-dimensional linear subspaces $L_v \oplus \CC\cdot v\subset \ten = \CC^{10|16}$
  are known as {\em super-null-geodesics} \cite{deligne-freed, kapranov-manin, movshev-schwarz, witten-10}. 
 \end{exas}

 It seems that 
 the $\PP(L_v)\cap\Sigma^{10}$ are precisely the { maximal complete intersection slices}
 of $\Sigma^{10}$. 
 
 Assuming this, we can formulate the constraints on superfields as follows.
 
 \begin{parab} {\bf Spin 0 constraints} are imposed on scalar superfields which are functions on $\Sc$
 or, more generally    maps $\Phi: \Sc\to X$ where $X$ is a given target manifold. They
 have the  form
 \[
 D_{b}\Phi = 0, \quad b\in B', 
 \]
 where $B'$ is one fixed  subspace of $B$ on which $\Gamma$ vanishes (case (0) above),
 maximal with this property. Maps $\Phi$ satisfying the constraints are known as {\em chiral superfields}.
  \end{parab}

  \begin{parab}\label{parab:spin1}
  {\bf Spin 1 constraints} are imposed on gauge superfields (connections  $\nabla$ in principal bundles on $\Sc$).
  They have the form of integrability $ (F_\nabla)|_{g\cdot \hen} =  0$, where $\hen=\ten_{\Gamma'}$ runs over all null-subalgebras in
  $\ten_\Gamma$ and $g\cdot\hen$ is the left translation of $\hen$ in $\Sc$. In other words, 
  \[
  [\nabla_{b_1}, \nabla_{b_2}] = \nabla_{\Gamma'(b_1, b_2)}, \quad [\nabla_b, \nabla_a] = [\nabla_{a_1}, \nabla_{a_2}]=  0, \quad
  \forall b,  b_1, b_2\in B', \,\,\, a, a_1, a_2\in V'
  \]
  where $B'$ runs over all the maximal complete intersection slices of $Y$ (case (1) above). 
  \end{parab}
  
   \begin{rem}
   When imposing constraints on superfields, it is obviously desirable not to end up restricting
   their dependence in the usual, even directions of spacetime. In the case  \ref{parab:spin1}
   this is ensured by the fact that the null-subalgebras $\hen$ have $\dim(\hen^\ev)=1$
   (they are the usual null-lines). 
   In other words,  {\em all the complete intersection slices of $Z=\Sigma^{10} $ 
    are quadric hypersurfaces}: $Z\cap M$ is a hypersurface in $M$.  
 It would be interesting to study other  intersections of quadrics  $Z$  with this property.  
 \end{rem}
 
 \begin{parab}
 {\bf Lie algebra meaning of complete intersections.}   
  Given a quadratic space $\Gamma: \Sym^2(B)\to V$, we can associate to it another,
 $\ZZ$-graded Lie algebra (i.e., a Lie algebra in $\Vect_\k^\ZZ$)
 \[
 \wt\ten_\Gamma \,\,=\,\, \on{FL}(B[-1]) \bigl/ (\Ker(\Gamma))
  \]
  Here $\FL(-)$ means the free graded Lie algebra generated by a graded vector space. 
  In our case $B[-1]$ is $B$ put in degree $+1$, so the degree 2 part of 
    $\on{FL}(B[-1])$   is
 $\Sym^2(B)$. The Lie algebra $\wt\ten_\Gamma$ is obtained by quotienting $\on{FL}(B[-1])$ by the graded Lie
 ideal generated by $\Ker(\Gamma)\subset\Sym^2(B)$. Thus
 \[
 \wt\ten^1_\Gamma = B, \,\,\, \wt\ten^2_\Gamma = V, 
 \]
  but it is not required that $V$ commutes with the generators and therefore $\wt\ten_\Gamma$ can be
  non-trivial in degrees $\geq 3$.  Cf. \cite{deligne-freed} \S 11.3. We
  note that  $\ten_\Gamma$ can be considered as a $\ZZ$-graded Lie algebra by lifting the degree $\ev$ part to degree $2$ and degree
  $\od$ part to degree $1$ (this is possible is possible since the degree $\ev$ part lies in the center). With this understanding, we
  have a surjective homomorphism
  of graded Lie algebras
  \[
p:   \wt\ten_\Gamma \lra \ten_\Gamma. 
  \]
  Denote by
  \[
  R_\Gamma \,\, =\,\,\Sym^\bullet(B^*) \bigl/ I_2(Y_\Gamma)
  \]
  the  graded coordinate algebra (commutative in the usual sense) of the homogeneous intersection of quadrics $Y_\Gamma\subset B$.
  Then, the enveloping algebra $U(\wt\ten_\Gamma)$ is identified with $R_\Gamma^!$, the quadratic dual of the quadratic algebra $R_\Gamma$,
  see \cite{PP} for background. In particular, we have a  homomorphism of graded algebras
  \[
 \eta:  U(\wt\ten_\Gamma) \lra \Ext^\bullet_{R_\Gamma}(\k, \k). 
  \]
 The algebra $R_\Gamma$ is called {\em Koszul}, if  $\eta$ is an isomorphism. This is the case in all spinorial examples.
 The role of complete intersections from this point of view is as follows.
 \end{parab}
 
 \begin{prop}
 The following are equivalent, and if they are true,  then $R_\Gamma$ is Koszul: 
 \begin{enumerate}
 \item[(i)] $\Gamma$ is of complete intersection type, i.e., $Y_\Gamma$ is a complete intersection of quadrics.
 
 \item[(ii)] We have $\wt\ten_\Gamma^i=0$ for $i\geq 3$, i.e., the morphism  $p: \wt\ten_\Gamma\to \ten_\Gamma$
 is an isomorphism. In particular,
 the condition of commutativity of $\wt\ten_\Gamma^2 = V$ with $\wt\ten_\Gamma^1=B$ already  follows from the defining
 relations of $\wt\ten_\Gamma$. 
 \end{enumerate}

 \end{prop}
  
 \noindent {\sl Proof:} This is a particular case of the general principle in commutative algebra that  locally complete intersections are characterized by the cotangent complex being quasi-isomorphic  to a 2-term complex.  The special case of intersections of
 quadrics was studied in
  \cite{kapranov-quadrics}. \qed
  
  \vskip .2cm
  
  Further, in many cases (including those related to spinors) the algebra $\wt \ten_\Gamma$ can be identified with the amalgamated
  free product of all its null-subalgebras $\wt \ten_{\Gamma'}=\ten_{\Gamma'}$. This relates the
  integrability conditions on null-subalgebras with the Koszul duality point of view on constraints for SYM  advocated in
  \cite{movshev-schwarz}.

 \section{Homotopy-theoretic underpinnings of supergeometry}
 
 \subsection{The skeleton of the Koszul sign rule}
 To understand the nature of the Koszul sign rule \ref{para:koszul}, let us ``minimize'' the symmetric monoidal category
 $\SVect_\k$ incorporating it. 
 
 To account just for the signs, we can disregard all morphisms in $\SVect_\k$ which are not isomorphisms,
 as well as all objects which have total dimension $>1$.  Restricted to
 1-dimensional super-vector spaces  ($\dim(V^\ev)+\dim(V^\od) =1$) and their isomorphisms, we get a symmetric monoidal category
 $\ISvect_\k$. Similarly, to capture the $\ZZ$-graded sign rule, we can restrict to the category $\IVect_\k^\ZZ$
 of $\ZZ$-graded 1-dimensional vector spaces. These categories are examples of the following concept.

 \begin{Defi}
 A {\em Picard groupoid} is a symmetric monoidal category $(\Gc, \otimes, \1, R)$ in which all objects
 are invertible  under  $\otimes$ and all morphisms are invertible  under composition. 
  \end{Defi}
  
  \noindent A Picard groupoid $\Gc$ gives rise to two abelian groups:
  
  \begin{itemize}
  \item The {\em Picard group} of $\Gc$, denoted $\on{Pic}(\Gc)$, or $\pi_0(\Gc)$. It is formed by
  isomorphism classes of objects, with the operation given by $\otimes$. 
  
  \item The group $\pi_1(\Gc) = \Aut_\Gc(\1)$ of automorphisms of the unit object. It is canonically identified
  with the group of automorphisms of any other object. 
  \end{itemize}

 In our case
 \[
 \beg
 \ISVect_\k: \quad \pi_0 = \ZZ/2, \quad \pi_1 = \k^*; 
 \\
 \IVect_\k^\ZZ: \quad \pi_0=\ZZ, \quad \pi_1=\k^*. 
 \eng
 \]
 Here, $\k^*$  is  still unnecessarily big: to formulate the sign rule, we need only the subgroup $\{\pm 1\}\subset\k^*$.
 So we   cut these Picard groupoids further. 
 
 For this, we replace $\k$ with the ring $\ZZ$, since
 $\{\pm 1\}=\ZZ^*$ is precisely its group of invertible elements. Accordingly,
  we replace 1-dimensional $\k$-vector
 spaces with free abelian groups of rank 1. This gives Picard groupoids $\ISAb, \IAb^\ZZ$. 
 Their objects are $\ZZ/2$-  or $\ZZ$-graded abelian groups which are free of rank 1. As before, the morphisms
 are isomorphisms, $\otimes$ is the graded tensor product over $\ZZ$ and the symmetries are given by
 the Koszul sign rule. The $\pi_i$ of these Picard groupoids are now as follows:
 \[
 \beg
 \ISAb: \quad \pi_0=\ZZ/2, \quad \pi_1=\{\pm 1\} = \ZZ/2, 
 \\
 \IAb^\ZZ: \quad \pi_0 = \ZZ, \quad \pi_1 = \{\pm 1\}=\ZZ/2. 
 \eng
 \]
 We can call $\ISAb$ and $\IAb^\ZZ$ the {\em sign skeleta} of the Koszul sign rule ($\ZZ/2$-graded and $\ZZ$-graded
 versions). They contain all the data needed to write the sign rule but nothing more. 
 
 The following simple but remarkable fact can be seen as a mathematical explanation of the  Principle
 of Naturality of Supers  \ref{para:nat-super}. 
 
 \begin{prop}\label{prop:freepic}
 $\IAb^\ZZ$ is equivalent to $\Fc_L$, the free Picard groupoid generated by one formal object (symbol)  $L$. 
 \end{prop}
 
 By definition, $\Fc_L$ has, as objects, formal tensor powers $L^{\otimes n}, n\in \ZZ$.
  It further has only those morphisms that are needed to write the symmetry isomorphisms
 \[
 R_{L^{\otimes m}, L^{\otimes n}}: L^{\otimes m}\otimes L^{\otimes n}= L^{\otimes (m+n)}
  \lra L^{\otimes n}\otimes L^{\otimes m} = L^{\otimes(n+m)}
 \]
 satisfying the axioms of a symmetric monoidal category (as well as composition, tensor products etc. of such
 morphisms). 
 
 \vskip .2cm
 
 \noindent  {\sl Sketch of proof:}  $L$ corresponds to the group $\ZZ$ placed in degree $1$, so $L^{\otimes n}$ is
 $\ZZ$ in degree $n$. Further, $R_{L,L}\in\Aut(L^{\otimes 2})=\Aut(\1)$ corresponds to $(-1)\in\ZZ^*$
 (note that $R_{L,l}\circ R_{L,L}=\Id$ by symmetry. The axioms of a symmetric monoidal category give that
 $R_{L^{\otimes m}, L^{\otimes n}}$ corresponds to $(-1)^{mn}$, so we recover the Koszul rule. \qed
 
 \vskip .2cm
 
 In other words, the category $\Vect_\k^\ZZ$ which is at the basis of all supergeometry,  can be obtained
 as a kind of  {\em $\k$-linear envelope of a free Picard groupoid}. More precisely, we have the following construction.
 
 \begin{Defi}\label{def:gradpic}
 Let $\Gc$ be a Picard groupoid, and $\chi: \pi_1(\Gc)\to \k^*$ be a homomorphism. By a {\em $(\Gc, \chi)$-graded
 $\k$-vector space} we will mean a functor $V: \Gc\to \Vect_\k$, whose value on objects
 will be denoted $A\mapsto V^A$, satisfying the following condition. For each object $A$, the action of each
 $\lambda\in\pi_1(\Gc) \simeq\Aut(A)$ on $A$ is taken into the multiplication by $\chi(\lambda)$ on $V^A$. 
 We denote by $\Vect_\k^{(\Gc,\chi)}$ the category of $(\Gc,\chi)$-graded $\k$-vector spaces. 
 \end{Defi} 
 
 Since $V$ is a functor, the spaces $V^A$ and $V^{A'}$ for isomorphic objects $A, A'$ are identified,
 so a $(\Gc, \chi)$-graded
 $\k$-vector space $V$ can be viewed as a $\pi_0(\Gc)$-graded vector space in the usual sense. 
 
 \begin{prop}\label{prop:gradpic}
 (a) The category $\Vect_\k^{(\Gc,\chi)}$ has a structure of a   monoidal category with the operation given
 by
 \[
 (V\otimes W)^A \,\,=\,\,\varinjlim\nolimits^{\Vect_\k}_{\{ B\otimes C\to A\}}  V^B \otimes_\k W^C, 
 \]
 the colimit taken over the category formed by pairs of objects $B, C\in \Gc$ together with an (iso)morphism
 $B\otimes C\to A$. Further, the symmetry in $\Gc$ makes $\Vect_\k^{(\Gc,\chi)}$ into a symmetric
 monoidal category. 
 
 \vskip .2cm
 
 (b) If $\Gc=\Fc_L$ and $\chi: \pi_1(\Fc_L)=\ZZ/2\to \k^*$ is the embedding of $\{\pm 1\}$, then
 $\Vect_\k^{(\Gc,\chi)}$ is identified with the category $\Vect_\k^\ZZ$ with the symmetry given by the Koszul sign rule. 
 \end{prop}
 
 \vskip .2cm
 
 
 \subsection{ (Higher) Picard groupoids and spectra}
 
 One of the insights of Grothendieck in  his manuscript ``Pursuing stacks"  (cf.\cite{deligne-det}, p. 114)
 was the correspondence between
 Picard groupoids and a particular class of {\em spectra} in the sense of homotopy topology.
 See \cite{Drinfeld} for a discussion and \cite{GK} for a slightly more detailed treatment which we follow here.
 A  systematic account can be found in \cite{JO}.
 
 \vskip .2cm
 
 The concept of a spectrum arises as a result of stabilizing the homotopy category of pointed topological spaces
 (say CW-complexes) under the two operations (adjoint functors)
 
 \[
 \beg
 \Sigma = \text{reduced suspension}, \quad \Omega = \text{loop space},
 \\
 \Hom(\Sigma X, Y) \,\,=\,\,\Hom(X, \Omega Y). 
 \eng
 \] 
 For example,  the spheres $S^n$ satisfy $\Sigma(S^n)=S^{n+1}$. We always have a canonical map
 (unit of adjunction)
 \[
 \eps_X: X\lra\Omega\Sigma X. 
 \]
 A spectrum $Y$ can be seen as a topological space $\Omega^\infty Y$ together with a sequence of
 deloopings: spaces $\Omega^{\infty-j}Y$ equipped with compatible homotopy equivalences
 $\Omega^j(\Omega^{\infty-j}Y)\sim\Omega^\infty Y$. A spectrum $Y$ has homotopy groups
 $\pi_i(Y), i\in \ZZ$ defined by
 \[
 \pi_i(Y) \,\,=\,\, \pi_{i+j}(\Omega^{\infty-j}Y), \quad j\gg 0. 
 \]
 
 \begin{ex}
 A topological space $X$ gives the {\em suspension spectrum} $\Sigma^\infty X$, with
 $\Omega^\infty\Sigma^\infty X = \varinjlim_n \Omega^n\Sigma^n X$
 and $\Omega^{\infty-j}\Sigma^\infty X=\varinjlim_n \Omega^{n-j}\Sigma^n X$  (limits under powers of $\epsilon$).
 The homotopy groups of $\Sigma^\infty X$ are the {\em stable homotopy groups} of $X$:
 \[
 \pi_i(\Sigma^\infty X) \,\,=\,\,\pi_i^\st (X) \,\, := \,\, \varinjlim_n \, \,  \pi_{i+n}\Sigma^n X. 
 \]
 \end{ex}
 
 Spectra form (after inverting homotopy equivalences), a triangulated category $\SHo$ known
 as the {\em stable homotopy category}. This category has a symmetric monoidal structure
 (smash product of spectra). 
 Let $m\leq n$ be integers ($m=-\infty$ or $n=\infty$ allowed).
 By an  $[m,n]$-{\em spectrum} we mean a spectrum $Y$ with
 $\pi_i(Y)=0$ with $i\notin [m,n]$, and  we denote by $\SHo_{[m,n]}\subset\SHo$ the full subcategory
 of $[m,n]$-spectra. There is a canonical ``truncation" functor
 \[
 \tau_{[m,n]}: \SHo\lra\SHo_{[m,n]}. 
 \]
 
 Grothendieck's correspondence can  be formulated as follows.
 
\begin{thm}
 There is an equivalence of categories
 \[
 \BB: \bigl\{ \text{Picard groupoids}\bigr\} [\on{eq}^{-1}] \lra \SHo_{[0,1]},
 \]
 so that $\pi_i(\BB(Gc) = \pi_i(\Gc)$, $i=0,1$. Here $ [\on{eq}^{-1}]$ means that equivalences of Picard groupoids
 are inverted, similarly to inverting homotopy equivalences in forming $\SHo$. 
 \end{thm} 
 
 A more precise result is proved in \cite{JO}. The spectrum $\BB\Gc$ corresponding to a Picard groupoid
 $\Gc$, is a version of the classifying space of $\Gc$. That is, the space $\Omega^\infty\BB\Gc=B\Gc$
 is the usual classifying space of $\Gc$ as a category, and the deloopings are constructed using the
 symmetric monoidal structure, see \cite{GK} (3.1.6) for an explicit construction. 
 
 \vskip .2cm
 
 The further point of Grothendieck is that more general $[0,n]$-spectra should have a description
 in terms of {\em Picard $n$-groupoids}, an algebraic concept to be defined, meaning
 ``symmetric monoidal $n$-categories with all the objects and higher morphisms invertible in all
 possible senses''.   Here we can  formally allow the case $n=\infty$.

 Incredible complexity of the stable homotopy category, well known to topologists, prevents us from
 hoping for a simple  algebraic definition of Picard $n$-groupoids. Nevertheless, for small values of
 $n=2,3$, this can be accessible and useful. The case $n=2$ is being treated in  the paper \cite{GJO}
 building on the theory 
  of symmetric monoidal 2-categories \cite{GJO-K, gurski-osorno}.

 
 \subsection {The sphere spectrum and the free Picard $n$-groupoid}

The fundamental role in homotopy theory is played by the {\em sphere spectrum} $\SSS = \Sigma^\infty S^0$
defined as the suspension spectrum of the $0$-sphere.  Its homotopy groups are the {\em stable homotopy
groups of spheres} 
\[
\pi_i(\SSS) \,\,=\,\,\pi_i^\st\,\,:=\,\, \varinjlim_n \pi_{i+n}(S^n), 
\]
 which vanish for $i<0$, so   $\SSS$ is a $[0,\infty]$-spectrum.  The spectrum $\SSS$
  is the unit object in the symmetric monoidal structure on 
 $\SHo$ and for this
reason can be considered as a homotopy-theoretic analog of the ring $\ZZ$ of integers.

 This motivates the further installment in Grothendieck's vision
 of a dictionary between spectra and Picard $n$-groupoids: 
   
 \begin{conj}\label{conj:freepic}
 The Picard
 $n$-groupoid corresponding to $\tau_{[0,n]}\SSS$, should be identified with $\Fc_L^{(n)}$ the free Picard $n$-groupoid
 on one formal object  $L$. 
 \end{conj}
 
 The concept of a free Picard $n$-groupoid presumes that we already have a system of axioms for
 what a Picard $n$-groupoid is. If we have such axioms, then $\Fc_L^{(n)}$ contains as objects,
 formal tensor powers $L^{\otimes n}$ and only those higher morphisms which are needed to write
 the necessary ``higher symmetry isomorphisms". 
 
 As before, for large $n$ this seems unattainable directly, but 
 for small $n\leq 2$ this can be made into a theorem. 
 In particular, the case $n=1$, proved in \cite[Prop. 3.1]{JO}  has enormous significance
 for super-geometry. Indeed, combining it with Proposition \ref {prop:freepic}, we arrive at:
 
 \begin{cor}
 The sign skeleton $\ISAb^\ZZ$ of the $\ZZ$-graded Koszul sign rule, is the Picard groupoid corresponding
 to the $[0,1]$-truncation of $\SSS$. The groups $\ZZ=\pi_0(\IAb^\ZZ)$ and $\ZZ/2=\pi_1(\IAb^\ZZ)$
 are the first two stable homotopy groups of spheres. 
 \end{cor}
 
 In other words, {\em the entire super-mathematics is obtained by unraveling  the first two layers of the
 sphere spectrum}.

 \vskip .2cm
 
 In Table \ref{table:corr} (which expands,  somewhat, a table from  the  online encyclopedia
 {\tt nLab})
 we give the values of the $\pi_i^\st$ for $i\leq 3$ and indicate mathematical and
 physical phenomena that these groups govern.
 We also compare $\SSS$ with another spectrum, the {\em algebraic K-theory spectrum}
 $\KK(\k)$ of a field $\k$, which has $\pi_i(\KK(\k))=K_i(\k)$,  the Quillen $K$-groups. These groups are indicated
 at in the bottom row. We notice that the first two groups are $\pi_i$ of our intermediate
 Picard groupoid $\IVect_\k^\ZZ$ which corresponds to the $[0,1]$-truncation of $\KK(\k)$, see
 \cite{deligne-det}, \S 4.
 A philosophy going back to \cite{BP} and to Quillen,  says that $\SSS$ can be  heuristically considered as the K-theory spectrum of  $\FF_1$, the
 (non-existent) field with one element, the symmetric group $S_n$ being the  ``limit'', as $q\to 1$, of the general
 linear groups $GL_n(\FF_q)$.

    \begin{table}
 {\scriptsize
  \renewcommand{\arraystretch}{2}
  \begin{tabular}{ |p{3cm} | p{3cm}| p{3cm} | p{3cm} |}
  \hline
 $\pi_0^\st=\ZZ$ & $\pi_1^\st=\ZZ/2$ & $\pi_2^\st=\ZZ/2$ & $\pi_3^\st=\ZZ/24$
 \\
 \hline
 Cardinality of sets, 
  dimension  of vector spaces&  Sign of permutation, determinant, Koszul rule & Spin, central extensions of symmetric groups
   & ``String'', $\sqrt[24]{1}$ in transformation of $\eta(q)$,
  $\chi(\on{K3})=24$
  \\
  \hline
  $K_0(\k)=\ZZ$ & $K_1(\k)=\k^*$ & $K_2(\RR)\to \pi_1(SL_n(\RR))$
  $\pi_1(SL_n(\RR))= \ZZ/2$ & Dilogarithms etc. 
  \\
  \hline
   \end{tabular}
       \caption{ 
       The first few $\pi_i^\st$ and their significance.}\label{table:corr}
}\end{table}

\vskip .2cm

The first two columns are self-explanatory.
In the third column,  the phenomenon  of { spin(ors)} is based on the fundamental group
$\pi_1(SO_n) = \pi_1(SL_n(\RR)) = \ZZ/2$ which is the same as $\pi_1^\st (SO_n)$
and is identified with  $\pi_1^\st=\ZZ/2$ via the map $SO_n\to\Omega^{n+1}(S^{n+1})$
(this is known as J-homomorphism).  
The existence of central extensions of symmetric and alternating groups $S_n$ and $A_n$
(with center $\ZZ/2$)
and of corresponding projective representations \cite{kleshchev, schur}  is a related phenomenon: $A_n$ embeds into $SO_n$,
and taking the preimage in $\Spin(n)$, we get a $\ZZ/2$-extension. 

One thing is worth noticing. 
 Supergeometry, as understood by mathematicians, tackles only the
first two columns of Table \ref{table:corr}.  A similar-sounding concept 
(supersymmetry) used by physicists,  dips
into the third column as well: fermions are always wedded to spinors in virtue of the Spin-Statistics Theorem. 
 In fact, there is something in the very structure of the sphere spectrum that seems to relate spin (third column) and
 statistics (the second column). At the most naive level, this is the coincidence of   $\pi_1^\st=\ZZ/2$
 with $\pi_2^\st=\ZZ/2$. 
 
 Further, consider the $[1,2]$-truncated spectrum $\tau_{[1,2]}\SSS$. Its loop spectrum
 $\Omega \tau_{[1,2]}\SSS$ is a $[0,1]$-spectrum with $\pi_0=\pi_1^\st=\ZZ/2$ and $\pi_1=\pi_2^\st =\ZZ/2$.

 \begin{para}\label{para:HSST}
 {\bf Homotopy-theoretic Spin-Statistics Theorem.}
 The Picard groupoid corresponding to $\Omega \tau_{[1,2]}\SSS$, is equivalent to $\ISAb$, the sign skeleton
 of the $\ZZ/2$-graded Koszul sign rule. In other words, $\Omega\tau_{[1,2]}\SSS$ is homotopy equivalent
 to  $(\tau_{[0,1]}\SSS)/2$,
 the  reduction of the spectrum $\tau_{[0,1]}\SSS$ by the element $2$ of its $\pi_0$. 
 \end{para}
 
 So there is not one, but {\em two ways} in which the same Koszul sign rule appears out of the sphere spectrum,
 one through statistics, the other one through spin.  Note that  the ``topological proof'' of the  usual
 physical Spin-Statistics Theorem, going back to Feynman \cite{feynman}  is based on  the intuitive claim that interchanging two particles is ``equivalent''  to  tracing a non-trivial loop in the rotation group,
 and this claim needs something like \ref{para:HSST} to be consistent. 
 See   \cite{duck}
 Ch.  20 for a  more detailed discussion  of this claim and the whole issue.

 \vskip .2cm
 
\noindent {\sl Proof of \ref{para:HSST}:}  This is an  exercise on known facts in homotopy theory.
A $[0,1]$-spectrum $Y$  (or a Picard groupoid) is classified by its $\pi_0, \pi_1$ and the
{\em Postnikov invariant} which is a group homomorphism $k_Y: \pi_0(Y)\to\pi_1(Y)$ satisfying $2k_Y=0$.
Explicitly, $k_Y$ is given by  the composition product  with the generator $\eta\in\pi_1^\st$.
For $Y=\tau_{[0,1]}\SSS$ we have therefore that $k_Y: \ZZ\to\ZZ/2$ is the surjection,
and for $Y=\tau_{[0,1]}\SSS/2$ we have  $k_Y = \Id_{\ZZ/2}$.

If both $\pi_0=\pi_1=\ZZ/2$ and  $k_Y\neq 0$,
 then $k_Y=\Id$, like for $Y=(\tau_{[0,1]}\SSS)/2$ and so $Y\sim (\tau_{[0,1]}\SSS)/2$.
Now, for $Y=\Omega\tau_{[1,2]}\SSS$ the Postnikov invariant is the map $\pi_1^\st\to\pi_2^\st$
given by composition with $\eta$. It is known   \cite{whitehead} that  $\eta^2\in\pi_2^\st$ is the generator
and so $k_Y\neq 0$. 
 \qed
 
 \vskip .2cm
 
 The 4th column of Table 1, headed by $\pi_3^\st=\ZZ/24$, is related to various ``string-theoretic"
 mathematics such as the appearance of 24th roots of 1 in Dedekind's formula for the modular transformation
 of his $\eta$-function \cite{atiyah}, the Euler characteristic of a K3 surface being 24, 
 the importance of the central charge modulo 24 in conformal field theory and so on, cf. \cite{H}.

 
 \subsection{Towards higher supergeometry}\label{sec:higher-super}
 
 Conjecture \ref{conj:freepic} means that the sphere spectrum $\SSS$ is the ultimate source for meaningful
 twists of commutativity, i.e., for designing  truly commutative-like structures, flexible enough
 to serve as a basis of geometry. The existing super-mathematics uses only the first two levels ($0$th and $1$st)
 of $\SSS$, with physical applications exploiting the parallelism between the $1$st and the $2$nd levels.

 This opens up a fantastic possibility of {\em higher super-mathematics} which would use, as its
 ``sign skeleton", the spectrum $\SSS$
 in its entirety or, at least, the truncations $\tau_{[0,n]}\SSS$ and the  free Picard $n$-groupoids $\Fc^{(n)}_L$
 for as long as we can make sense of them algebraically. Here we sketch the first step in this direction,
 the formalism for  $n=2$.  For convenience, we adopt a genetic approach. 
 
 \begin{parab}\label{parab:idea}
 {\bf Idea of supersymmetric monoidal categories.}   Super-mathematics begins with replacing
 commutativity $ab=ba$ with supercommutativity \eqref{eq:supercomm}. 
 Categorical analogs of commutative algebras are
 symmetric monoidal categories, where we have coherent isomorphisms $V\otimes W \simeq W\otimes V$. 
 So we introduce  ``categorical minus signs" into these isomorphisms as well. 
 
More precisely, by a $\k$-{\em superlinear category} we mean a module category $\Vc$ (a category tensored over)
the symmetric monoidal category $\SVect_\k$. 
In such a category we have the {\em parity change functor} $\Pi$ given by tensoring with the super-vector space
$\k^{0|1}$. We take $\Pi$ as the categorical analog of the minus sign. In doing so, we use the identification
of $\pi_0(\SSS/2)=\ZZ/2$ (the $\ZZ/2$-grading) with $\pi_1^\st = \ZZ/2$ (the $\pm 1$ signs). 

We now consider $\k$-superlinear categories $\Ac$ which are  $\ZZ$- or $\ZZ/2$-{\em graded}, i.e., split into a
categorical direct sum $\Ac = \boxplus_i \Ac^i$, where $i\in \ZZ$ or $\ZZ/2$. We assume that $\Ac$
is equipped with graded $\SVect_\k$-bilinear bifunctors
\[
\otimes \, =\, \otimes_{i,j}: \Ac^i \times \Ac^j \lra\Ac^{i+j}
\]
subjects to associativity isomorphisms of the usual kind and want to impose {\em twisted commutativity isomorphisms}
\be\label{eq:twisted-com}
R_{V,W}: V\otimes W \lra \Pi^{\deg(V)\cdot\deg(W)}  (W\otimes V) 
\ee
subjects to natural axioms, in which, further,  various numerical minus signs will be introduced. 
  
 \end{parab}
 
 \begin{parab} {\bf Definition of supersymmetric monoidal categories.} 
 A possible more precise definition can go as follows. 
 For simplicity  consider  the $\ZZ$-graded version. 
 We use Proposition \ref{prop:gradpic} as a guideline, and start with $\Fc= \Fc^{(2)}_L$, the free Picard 2-groupoid
 on one object $L$. It corresponds to the truncation $\tau_{[0,2]}\SSS$. 
 Thus $\pi_0(\Fc)$, the group of equivalence classes of objects, is  identified with $\pi_0^\st=\ZZ$
 and will account for the grading. 
 
 The category $\Aut_{\Fc}(\1)$ formed by automorphisms of the unit object
 and 2-morphisms between them, is a usual Picard groupoid which corresponds to $\Omega\tau_{[1,2]}\SSS$
 and so, by the  "Spin-Statistics Theorem''
 \ref{para:HSST},  it is identified with $\ISAb$, the sign skeleton of $\SVect_\k$.
 
 We denote by $\SCat_\k$ the 2-category of $\k$-superlinear categories.
 It serves as a categorical analog of the category of ordinary vector spaces. We will denote by
 $\Vc\boxtimes_{\SVect-\k}\Wc$ the categorical tensor product \cite{greenough}
 of two $\k$-superlinear categories $\Vc$ and $\Wc$.
 
 The Picard groupoid
 $\ISVect_\k$ plays the role of the multiplicative group $\k^*$ for $\SCat_\k$: it acts on each object
 by equivalences. The monoidal functor (embedding) ${\boldsymbol \chi}: \ISAb\to\ISVect_\k$
 is  therefore an analog of the homomorphism $\chi$ from
  Definition \ref{def:gradpic}. 
  
  We now consider the 2-category $\SVect^{(\Fc, \chib)}_\k$
 formed by all 2-functors $\Vc: \Fc\to\SCat_\k$ which take the action of $\Aut_{\Fc}(\1)$
 on each object $A$ into the action on $\Vc^A$ given by $\chib$. 
 As before, a datum of such $\Vc$ is the same as a datum of a family of superlinear categories
 $\Vc^i = \Vc^{L^{\otimes i}}$, $i\in \ZZ$, one for each equivalence class of objects of $\Fc$. 
Now, the formula
\[
(\Vc \boxtimes \Wc)^A \,\,=\,\, 2\varinjlim\nolimits_{\{ B\otimes C\to A\}} \Vc^B \boxtimes_{\SVect_\k} \Wc^C
\]
makes $\SVect^{(\Fc, \chib)}_\k$ into a symmetric monoidal 2-category. It can be seen as the categorical
analog of the  category of super-vector spaces.

By definition, a {\em supersymmetric monoidal $\k$-category} is a symmetric monoidal  object $\Ac$ in 
$\SVect^{(\Fc, \chib)}_\k$.  An  explicit algebraic model for $\Fc$ was proposed 
in \cite[Ex. 5.2]{bartlett}, see also \cite[Ex. 2.30]{schommer-pries}.  Taking this model for $\Fc$, we can unravel the data
involved in $\Ac$.  These data in particular,  contain  superlinear categories $\Ac^i$,  bifunctors $\otimes_{i,j}$
and isomorphisms $R_{V,W}$ as  outlined in \ref{parab:idea}. 

The definition of a $\ZZ/2$-graded supersymmetric monoidal category is similar, using the Picard
2-groupoid $\Fc/2$. 

\vskip .2cm

 If $V$ is an object of a  symmetric monoidal category, then there is an action of the symmetric group $S_n$
 on $V^{\otimes n}$. If, instead $V\in\Ac^{2m+1}$ is an odd object of a supersymmetric monoidal category $\Ac$,
 then $V^{\otimes n} \oplus\Pi(V^{\otimes n})$
 has  an action of the (spin)  central extension of $S_n$, first discovered by Schur \cite{schur}. 
 
  \end{parab}
  
  \begin{exas} {\bf (a) The exterior algebra of a superlinear category.} Similarly to usual linear algebra (Examples \ref{exas:supercomm}),
  an example of a supersymmetric monoidal category can be extracted from the {\em categorical version of the exterior power
  construction} developed  in \cite{GK}. This construction is based on the {\em categorical sign character}
  which is a functor of monoidal categories
  \[
\on{sgn}_2: S_n\lra \ISAb
  \]
  ($S_n$ is the symmetric group considered as a discrete monoidal category). It combines the usual sign character
  $\on{sgn}: S_n\to \ZZ/2$ at the level of $\pi_0$ and the ``spin-cocycle" $c\in H^2(S_n, \ZZ/2)$ at the level of $\pi_1$.
  The exterior power $\Lambda^n\Vc$ of a superlinear category $\Vc$ is obtained from the  tensor power
  $\Vc^{\boxtimes n}$  by considering objects equipped with  $\on{sgn}_2$-twisted 
  $S_n$-equivariance structure (see \cite{GK} for precise context and details). The analog of the wedge product is given by
  the functors
  \[
  \wedge_{m,n}: \Lambda^m\Vc \times\Lambda^n\Vc \lra\Lambda^{m+n}\Vc
  \]
  given by partial $\Pi$-antisymmetrization, as in \cite{GK} \S 4.2. 
  
  \vskip .2cm
  
 {\bf  (b) Superalgebras of types M and Q and the half-tensor product of Sergeev.}
  Let $\k$ be algebraically closed. It is known since C.T.C. Wall \cite{wall} that simple  finite-dimensional
  associative superalgebras over $\k$ are of two types:
  {\em type M}, formed by the  matrix superalgebras $M_{p|q}=\End(\k^{p|q})$ and {\em type Q}, formed by the so-called
  {\em queer superalgebras} $Q_n\subset M_{n|n}$,  see \cite{jozefiak} and 
  \cite{kleshchev}. The simplest nontrivial queer algebra is the Clifford algebra $\on{Cliff}_1$ on one generator
  \[
  Q_1 = \on{Cliff}_1, \quad \on{Cliff}_n \,\,  = \,\,\k[\xi_1, \cdots\xi_n]/\bigl(\xi_i^2=1, \,\, \xi_i\xi_j=-\xi_j\xi_i\bigr), \,\,\,\deg(\xi_i)=\od. 
  \]
  Their behavior under tensor multiplication is
  \be\label{eq:MQ}
  \begin{gathered}
  M_{p|q}\otimes M_{m|n}\,\, \simeq\,\, M_{pm+qn|pn+qm}, 
  \\
 M_{p|q}\otimes Q_n
 \,\simeq\, Q_{(p+q)n}, \quad Q_m\otimes Q_n\, \simeq \, M_{mn| mn}.  
  \end{gathered}
  \ee
 This means that the {\em super-Brauer group} formed by  Morita equivalence classes of these algebras,  is identified with $\ZZ/2$,
 with type $M$ mapping to $\ev$ and type $Q$ mapping to $\od$.  
 This $\ZZ/2$ is nothing but $\pi_2^\st$, responsible for spin.  See \cite{Freed-anomalies} \S 4. 
 
 As a consequence,  irreducible objects of any semisimple  $\k$-superlinear category $\Vc$  also split  into  two types $M$ and $Q$,
 according to their endomorphism algebras being $\k$ or $Q_1$. Denoting $\Vc^\ev$ the subcategory formed by direct sum
 of objects of type $M$ and $\Vc^\od$ the subcategory formed by sums of objects of type $Q$, we get an intrinsic $\ZZ/2$-grading on $\Vc$. 
 By \eqref{eq:MQ}, any exact  monoidal structure $\otimes$ on $\Vc$ preserves this grading.
 
  Further, if $V, W$ are irreducible objects of
 type $Q$, then $V\otimes W$ is acted upon by $Q_1\otimes Q_1 \simeq \End(\k^{1|1})$ 
 and so is identified with the direct sum of some object with its shift:
 \[
V\otimes W \,\,\simeq \,\, ( 2^{-1}V\otimes W) \oplus \Pi(2^{-1}V\otimes W), \quad
2^{-1}V\otimes W :=   I \otimes_{Q_1\otimes Q_1} (V\otimes W) 
 \]
 where $I\simeq (\k^{1|1})^*$ is an irreducible right module over $Q_1\otimes Q_1$.
   We get in this way a new monoidal operation
 \[
 \otimes_{\od, \od}: \Vc^\od \times \Vc^\od \lra \Vc^\ev, \quad V\otimes_{\od, \od} W := 2^{-1}V\otimes W.
 \]
 This operation was  introduced by A. Sergeev \cite{sergeev-q}, see also \cite{kleshchev}, p.163  for more discussion. 
  If $\otimes$ is symmetric, then $\otimes_{\od, \od}$ satisfies
 \[
 V\otimes_{\od,\od} W \,\,\simeq \,\,  \Pi(W\otimes_{\od, \od}V).
 \]
 Indeed, the  interchange of the factors in $V\otimes W$ corresponds to the  interchange of $\xi_1$ and $\xi_2$ in $\on{Cliff}_2=Q_1\otimes
 Q_1$, and the pullback of $I$ under this interchange is isomorphic to $\Pi(I)$ (which, as a right $\on{Cliff}_2$-module, is not
 isomorphic to $I$). 
  The operation  $\otimes_{\od, \od}$  can be used as a source of examples of supersymmetric monoidal structures. 
  \end{exas}

  \begin{parab} {\bf Further directions.} It seems important to investigate the concept of
  supersymmetric monoidal categories more systematically.  Among other things, in such a  category $\Vc$
  we can speak about commutative algebra objects $A$. They would be further  twisted generalizations
  of super-commutative algebras (which appear for $\Vc=\SVect_\k$,   a symmetric, not just
 a  supersymmetric category). At the same time, such objects should be ``just as good''
 as  ordinary commutative algebras, in an extension of the Principle of Naturality of Supers
  \ref{para:nat-super}. In particular, we must be able to associate to them some geometric
  objects $\Spec(A)$ which should be flexible enough to be glued together to  form more global
  structures. 
  
  \vskip .2cm
  
  It is also interesting   to weaken the concept to allow, e.g.,  for  ``super-braidings'' etc.
  in which the twisted commutativity isomorphisms   \eqref{eq:twisted-com} 
  are retained but the involutivity requirements on them are relaxed. 
  
  \vskip .2cm 
  
 Super-symmetric and super-braided
 monoidal categories may be relevant for the
  recent activity in condensed matter physics studying so-called
 fermionic phases, see \cite{gu-wen} \cite{Freed-entangl} \cite{Freed-Hopkins} \cite{bhardwaj}. 
 In particular, various ``super-cohomooogy" constructions appearing there, seem to be
 represented by homotopy types and spectra related to truncations of $\SSS$, the sphere spectrum and
 $\Kc(\CC)$, the K-hteory spectrum of $\CC$. For instance,  ``group super-cohomology''' of 
 \cite{gu-wen} is represented by the spectrum $\tau_{\leq 1}\Kc(\CC)/2$, which corresponds to
 the Picard groupoid of 1-dimensional super-vector spaces over $\CC$. 
 
 \vskip .2cm
 
 Another possible approach to higher super-geometry could be by using the truncations and suspensions not of $\SSS$  itselft but of
 \[
 \widecheck\SSS \,\,=\,\,  \ul{R\Hom} (\SSS, \CC^*),
 \]
 the Pontryagin (Brown-Comenetz) dual of $\SSS$ considered in  \cite{Freed-entangl}. The homotopy groups of $\widecheck\SSS$
 are nontrivial only in degrees $\leq 0$ with $\pi_0=\CC^*$ dual to $\pi_0^\st=\ZZ$ and other $\pi_i$ finite.  
 For example, the spectrum corresponding to the Picard groupoid of 1-dimensional super-vector spaces above can   be  also described  as
 \[
 \tau_{\leq 1} \Kc(\CC)/2 \,\,=\,\, \Sigma \tau_{[-1,0]}  \widecheck\SSS,.
 \]
 the suspension of the truncation of $\widecheck\SSS$. 
 
 \vskip .2cm
 
 As for the next step,  one can imagine a ``2-supersymmetric monoidal 2-category''  to split into a direct sum of  $24$  sub(2-)categories ( ``sectors'' )
     instead
  of just two   (even and odd objects), with    Hom-categories between  objects split into even and odd parts and so on,
  with an appropriate pattern of sign/shift  rules for symmetries.  The pattern of 24 sectors is  of course
   suggestive of the features of conformal field theory and
  Chern-Simons theory (the significance of the central charge modulo 24, the 24 classes of framings of a 3-manifold
  etc.).

 \vskip .2cm
 
 Expressed more generally,  study of higher supergeometry can include two complementary directions. 
 One is {\em categorification}, investigation of deeper and deeper twists for ``commutativity'',
 taking the inspiration from the ultimate commutative structure, the sphere spectrum $\SSS$.
 The other is  {\em geometrization}, gluing together local objects of algebraic origin, for which
 some sort of ``commutativity'' seems to be necessary (why?). 
 Reconcilng these two directions can, hopefully, shed more light on the way geometry
  emerges from quantum behavior.

  \end{parab}

   \vfill\eject

\let\thefootnote\relax\footnote {
 Kavli Institute for Physics and Mathematics of the Universe (WPI), 5-1-5 Kashiwanoha, Kashiwa-shi, Chiba, 277-8583, Japan.
Email: {mikhail.kapranov@ipmu.jp}

\vskip .2cm
}


\begin{thebibliography}{WWW}
 
 \small

 \bibitem{Ad} J. F. Adams. Infinite Loop Spaces.
 Princeton Univ. Press, 1978. 
 
 \bibitem{atiyah} M. Atiyah. The logarithm of the Dedekind $\eta$-function.
 {\em Math. Ann.} {\bf 278} (1987) 335-380. 
 
 \bibitem{ABS} M. Atiyah, R. Bott, A. Shapiro. Clifford modules.
 {\em Topology} {\bf 3} (1964) 3-38. 
 
 \bibitem{BP}
M. Barratt, S. Priddy. On the homology of non-connected monoids and their associated groups.
{\em  Comment. Math. Helv.}  {\bf 47} (1972), 1-14.

\bibitem{bartlett} B. Bartlett. Quasistrict symmetric monoidal 2-categories via wire 
diagrams. arXiv:1409.2148. 

\bibitem{bhardwaj} L.Bhardwaj, D. Gaiotto, A. Kapustin. State sum constructions of spin-TFT
and string net consntructions of fermionic phases of matter. arXiv:1605.01640. 

\bibitem{CFK} I. Ciocan-Fontanine, M. Kapranov. Derived $Quot$ schemes.
{\em Ann. Sci. ENS}, {\bf 34} (2001) 403-440. 

\bibitem{deligne-det} P.  Deligne.
Le d\'eterminant de la cohomologie, in: {\em Current trends in arithmetical algebraic geometry} (Arcata, Calif., 1985), 93-177,
Contemp. Math., {\bf 67} , Amer. Math. Soc., Providence, RI, 1987.

\bibitem{deligne-tannakiennes} P. Deligne. Cat\'egories Tannakiennes,
 Grothendieck Festschrift, vol. II, pp.111-195, Birkh\"auser, 1990. 
 
 

\bibitem{deligne-spinors} P. Deligne. Notes on spinors, in:  ``Quantum Fields and Strings: a Course for Mathematicians" 
(P. Deligne et al. Eds.) vol. I, p. 99-136, Amer. Math. Soc., Providence, RI, 1999.

\bibitem{deligne-tensor} P. Deligne.  Cat\'egories Tensorielles. {\em Moscow Math. J.,}
{\bf 2}  (2002)  227-248. 

\bibitem{deligne-morgan} P. Deligne, J. Morgan. Notes on supersymmetry  (following J. Bernstein), in:
  ``Quantum Fields and Strings: a Course for Mathematicians'' 
(P. Deligne et al. Eds.) vol. I, p. 41-98, Amer. Math. Soc., Providence, RI, 1999.

\bibitem{deligne-freed} P. Deligne, D. S. Freed. Supersolutions, in:
``Quantum Fields and Strings: a Course for Mathematicians" 
(P. Deligne et al. Eds.) vol. I, p.  227-356, Amer. Math. Soc., Providence, RI, 1999.

\bibitem{DW1}	R. Donagi, E. Witten.  Supermoduli space is not projected.
arXiv:1304.7798.
 
\bibitem{DW2} R. Donagi, E. Witten. 
 Super Atiyah classes and obstructions to splitting of supermoduli space.  	arXiv:1404.6257. 

\bibitem{Drinfeld} V. Drinfeld. Infinite-dimensional objects in algebra and geometry, in:
``The Unity of Mathematics (In Honor of the Ninetieth Birthday of I.M. Gelfand)'' pp. 263-304, Birkh\"auser, 2006. 

 

\bibitem{duck} I. Duck, E. C. G. Sudarshan. Pauli and the Spin-Statistics Theorem. World Scientific, Singapore, 1997. 

 

\bibitem{feynman} R. P. Feynman, S. Weinberg. Elementary particles and the laws of physics. The 1986 Dirac Memorial Lectures.
Cambridge Univ. Press, 1999. 


\bibitem{Freed-super} D. S. 
Freed. Five Lectures on Supersymmetry. Amer. Math. Soc. , Providence, RI, 1999.

\bibitem{Freed-anomalies} D. S. Freed. Anomalies and invertible field theories. ArXiv:1404.7224.

\bibitem{Freed-entangl}  D. S. Freed.  Short-range entanglement and invertible field theories. ArXiv:1406.7278. 

\bibitem{Freed-Hopkins}  D. S. Freed, M. J. Hopkins. Reflection positivity and invertible topological phases.
ArXiv:1604.06527. 


\bibitem{GK}
N. Ganter, M. Kapranov. Symmetric and exterior powers of categories.
{\em Transform. Groups} {\bf 19} (2014)  57-103.

\bibitem {GGRS} S. J. Gates, Jr., M. T.  Grisaru, M. Rocek, and W.  Siegel. Superspace, or 1001 Lessons in Supersymmetry. 
B. Cummings Publ. Reading, MA, 1983. 

\bibitem{greenough} J. Greenough. Monoidal 2-structure of bimodule categories.
arXiv:0911.4979. 

\bibitem{gu-wen} Z.-C. Gu, X.-G. Wen. Symmetry-protected topological orders for interacting fermions:
Fermionic topological nonlinear $\sigma$ models and a special group supercohomology theory. 
  {\em Phys. Rev..  B}  {\bf 90}  (2014) 115141. 

\bibitem{GJO-K} N. Gurski, N. Johnson, A.M. Osorno.  K-theory for 2-categories. ArXiv:1503.07824. 

\bibitem{GJO} N. Gurski, N. Johnson, A.M. Osorno. 
Realizing
stable 2-types via Picard 2-categories. In preparation. 

\bibitem{gurski-osorno} N. Gurski, A. M. Osorno.
Infinite loop spaces and coherence for symmetric monoidal bicategories.
{\em Adv. Math.} {\bf 246} (2013) 1-32. 

 \bibitem{H} M. J. 
 Hopkins. Algebraic topology and modular forms.
 Proceedings of the International Congress of Mathematicians, Vol. I (Beijing, 2002), 291-317, Higher Ed. Press, Beijing, 2002.

\bibitem{HS} M. J. Hopkins, I. Singer. Quadratic functions in geometry, topology and M-theory.
{\em J. Diff. Geom.} {\bf 70} (2005) 329-452.

\bibitem{JO}  N. Johnson, A.  M. Osorno. 
Modeling stable one-types.
{\em Theory and Applications of Categories}, {\bf 26}  (2012)   520-537.

\bibitem{jozefiak} T. Jozefiak. Semisimple superalgebras, in:
 ``Algebra-: some current trends" (Varna, 1986), 96-113, {\em Lecture Notes in Math.} {\bf 1352},
  Springer, Berlin, 1988. 
  
  \bibitem{kapranov-quadrics} M. M. Kapranov. On the derived category and K-functor of coherent sheaves
  on intersections of quadrics. {\em Math. USSR Izv.} {\bf 32} (1989) 191-204. 
  
 \bibitem{kapranov-manin} M. M. Kapranov, Y. I. Manin. 
 The twistor transformation and algebraic-geometric constructions of solutions of the equations of field theory.
 {\em Russian Math. Surveys,}  {\bf 41}:5 (1986)  33-61. 	
 
 \bibitem{KV} M. Kapranov, E. Vasserot. Supersymmetry and the formal loop space. 
{\em  Adv. Math.} 
{\bf 227} (2011) 1078-1128.
 

\bibitem{kleshchev} A. Kleshchev. Linear and Projective Representations of Symmetric Groups.
Cambridge Univ. Press, 2005. 

\bibitem{kontsevich} [28] M. Kontsevic. Notes on deformation theory, course lecture notes, Berkeley,
1992.

\bibitem{leites} D. A. Leites. Spectra of graded-commutative rings. 
 {\em  Uspekhi. Mat. Nauk}  {\em 29}:3  (1974) 209-210. 
 
 \bibitem{lott} J. Lott. Twistor constraints in supergeometry. {\em Comm. Math. Phys.}
 {\bf 133} (1990) 563-615.
 
 \bibitem{lurie} J. Lurie. Derived algebraic geometry II-V. arxiv math/0702299, math/0703204, 0709.3091, 0905.0459.  

\bibitem{maclane} S. Mac Lane. Categories for a Working Mathematician. Springer-Verlag, 1971. 

\bibitem{M} Y. I. Manin. Gauge Fields and Complex Geometry. Springer-Verlag, 1997. 

\bibitem{movshev-schwarz} M. Movshev, A. Schwarz.
On maximally supersymmetric Yang-Mills theories.
{\em Nucl. Phys. B} {\bf 681} (2004), 324-350. 

\bibitem{nahm} W. Nahm.  Supersymmetries and their representations.
{\em  Nucl. Phys. B}, {\bf 135} (1978) 149-166.

\bibitem{PP}  A. Polishchuk, L. Positselski. Quadratic Algebras. Amer. Math. Soc., Providence, RI, 2005. 

\bibitem{rama}  N.  Ramachandran.
 Values of zeta functions at  $s=1/2$.  {\em Int Math Res Notices} {\bf 25}  (2005) 1519-1541. 

\bibitem{reid} M. Reid. The complete intersection of two or more quadrics. Thesis, Cambridge University, 1972. 

\bibitem{Sato-theta} M. Sato. Pseudo-differential equations and theta functions. Colloque International
 CNRS sur les \'equations aux D\'eriv\'ees Partielles Lin\'eaires (Univ. Paris-Sud, Orsay, 1972), pp. 286- 291.
 {\em Ast\'erisque}, {\bf 2-3},  Soc. Math. France, Paris, 1973.
 
\bibitem{SKK-theta}
M. Sato, M. Kashiwara, T.  Kawai.  Linear differential equations of infinite order and theta functions.
{\em  Adv. in Math.}  {\bf 47}  (1983)  300-325.  


\bibitem{schommer-pries}  C. Schommer-Pries. The classification of two-dimensional extended topological field theories.
arXiv: 1112.1000. 

\bibitem{schur} I. Schur.   \"Uber die Darstellung der symmetrischen und der alternierenden Gruppe durch gebrochene lineare Substitutionen.
{\em  J. Reine Angew. Math.}  {\bf 139} (1911) 155-250. 


 \bibitem{sergeev-q} A. N. Sergeev. The tensor algebra of the identity representation as a module over the
 Lie superalgebras $\mathfrak{Gl}(m,n)$ and $Q(n)$.
 {\em Math. USSR. Sbornik}, {\bf 51} (1985) 419-427. 

\bibitem{Sihn} H. X. Sinh. Gr-cat\'egories. Th\`ese, Universit\'e Paris-7, 1975. 


\bibitem{stolz-teichner} S. Stolz, P. Teichner. Supersymmetric field theories and generalized cohomology, in:
``Mathematical Foundations of Quantum Field Theory and Perturbative String Theory'', 
{\em Proc. of Symp. in Pure Math.} {\bf  83}  (2011)  279-340, Amer. Math. Soc. , Providence, RI, 2011. 

\bibitem{T} Y. Tachikawa. A pseudo-mathematical pseudo-review on $4d$  $N = 2$  supersymmetric quantum field theories. 
IPMU preprint, 2014.

\bibitem{TV} B. To\"en, G Vezzosi.  Homotopical algebraic geometry II: geometric stacks and applications. 
{\em Mem. Amer. Math. Soc.}  {\bf 193} (2008), no. 902. 

\bibitem{tyurin} A. N. Tyurin. On intersections of quadrics. {\em Russian Math. Surveys}, {\bf 30}:6 (1975), 51-105. 

\bibitem{wall} C. T.  C. Wall. Graded Brauer groups. 
{\em J. reine und angew. Math.} {\bf 213} (1964) 187-199.  

\bibitem{WB} J. Wess, J. Bagger. Supersymmetry and Supergravity. Princeton University Press, 1991. 

\bibitem {whitehead} G.W.  Whitehead. Recent Advances in Homotopy Theory.   Amer. Math. Soc., Providence, RI,  2007. 

\bibitem{witten-10} E. Witten. Twistor-like transform in ten dimensions. {\em Nucl. Phys. B} {\bf 266} (1986) 245-264. 

\bibitem{witten-morse} E. Witten, Supersymmetry and Morse theory,  {\em J. Differential Geom.}  {\bf 17} 
(1982) 661?692.

 \end{thebibliography}
 \end{document}